\documentclass[gdplain]{geradwp}

\PassOptionsToPackage{hyphens}{url}
\graphicspath{{graphics/}} 
\usepackage{fancyhdr,graphicx,amsmath,amssymb}

\usepackage{subfigure}
\usepackage{multirow}
\usepackage{mathtools}

\usepackage[ruled,vlined]{algorithm2e}
\usepackage{hyperref}
\usepackage{natbib}
 \bibpunct[, ]{(}{)}{,}{a}{}{,}%

\usepackage{changepage}
\usepackage{amsmath}
\usepackage[T1]{fontenc}

\hypersetup{colorlinks,%
citecolor={blue}, 
urlcolor={blue},
linkcolor={blue},
breaklinks={true}
}
\newcommand{\evspp}{EVSPP-OI}

\newcommand{\chr}{Ch.}
\newcommand{\dd}{Dr.}

\newcommand{\dc}{Wa.}

\newcommand{\roi}{ROI}
\newcommand{\rss}{RSS}

\newcommand{\RNum}[1]{\uppercase\expandafter{\romannumeral #1\relax}}
\newcolumntype{C}[1]{>{\centering\arraybackslash}p{#1}}

\usepackage[normalem]{ulem}

\GDtitre{Dynamic Routing for the Electric Vehicle Shortest Path Problem with Charging Station Occupancy Information}
\GDmois{Avril}{April}
\GDannee{2023}
\GDnumero{53}
\GDauteursCourts{M. Dastpak, F. Errico, O. Jabali, F. Malucelli}
\GDauteursCopyright{Dastpak, Errico, Jabali, Malucelli}

\begin{document}

\GDpageCouverture

\begin{GDpagetitre}

\begin{GDauthlist}
\GDauthitem{Mohsen Dastpak \ref{affil:ets}\GDrefsep\ref{affil:gerad}\GDrefsep\ref{affil:cirrelt}}
\GDauthitem{Fausto Errico \ref{affil:ets}\GDrefsep\ref{affil:gerad}\GDrefsep\ref{affil:cirrelt}}
\GDauthitem{Ola Jabali \ref{affil:polimi}}
\GDauthitem{Federico Malucelli \ref{affil:polimi}}
\end{GDauthlist}

\begin{GDaffillist}
\GDaffilitem{affil:ets}{Department de génie de la construction, École de technologie supérieure, Montr\'eal (Qc), Canada, H3C 1K3}
\GDaffilitem{affil:gerad}{GERAD, Montr\'eal (Qc), Canada, H3T 1J4}
\GDaffilitem{affil:cirrelt}{CIRRELT, Montr\'eal (QC), Canada, H3C 3J7}
\GDaffilitem{affil:polimi}{Dipartimento di Elettronica, Informazione e Bioingegneria, Politecnico di Milano, Piazza Leonardo da Vinci 32, Milano 20133, Italy}
\end{GDaffillist}

\begin{GDemaillist}
\GDemailitem{mohsen.dastpak.1@ens.etsmtl.ca}
\end{GDemaillist}

\end{GDpagetitre}

\GDabstracts

\begin{GDabstract}{Abstract}
We study the problem of an Electric Vehicle (EV) having to travel from an origin to a destination in the shortest amount of time. We focus on long-distance settings, where the shortest path between the origin and the destination has energy requirements exceeding the EV autonomy. The EV may charge its battery at public Charging Stations (CSs), which are subject to unknown arrivals of exogenous vehicles requiring uncertain charging times. Thus, the waiting times at CSs are uncertain. Similar to other contributions in the literature, we model CSs using appropriately defined queues, whose status is revealed upon the EV arrival. However, following recent technological advances, we also consider that the status of each CS is updated in real-time via binary Occupancy Indicator (OI) information signaling if a CS is busy or not. Therefore, we assume that the EV continuously receives OI updates on all CSs. At each update, we determine the sequence of CSs to visit along with associated charging quantities. We name the resulting problem as the Electric Vehicle Shortest Path Problem with charging station Occupancy Indicator information (\evspp{}). In this problem, we consider that the EV is allowed to partially charge its battery, and we model charging times via piecewise linear charging functions that depend on the CS technology.

We propose an MDP formulation for the \evspp{}, which aims at optimizing the EV routing and charging policy. To solve the problem, we develop a reoptimization algorithm that establishes the sequence of CS visits and charging amounts based on system updates. Specifically, we  propose a simulation-based approach to estimate the waiting time of the EV at a CS as a function of its arrival time. As the path to a CS may consist of multiple intermediate CS stops, estimating the arrival times at each CS is fairly intricate. 
To this end, we propose an efficient heuristic that yields approximate lower bounds on the arrival time of the EV at each CS, which are used to derive an estimation of the waiting time at each CS. We use these estimations to define  a compatible  deterministic version of the EVSPP, which we solve with an existing algorithm. 
We conduct a comprehensive computational study and compare the performance of our methodology with a benchmark that observes the status of CSs only upon arrival (i.e., with no OI information). Results show that our method reduces waiting times and total trip duration by an average of 23.7\%-95.4\% and 1.4\%-18.5\%, respectively.

\end{GDabstract}

\GDarticlestart

\section{Introduction}
Electric Vehicles (EVs) have taken center stage in recent years and have the potential to revolutionize the way we travel and interact with our environment.
BloombergNEF (Bloomberg New Energy Finance) projects that plug-in vehicle sales will rise from 6.6 million in 2021 to 20.6 million in 2025~\citep{BloombergNEF2022}.
While this trend will lead to more sustainable transportation, a number of technical barriers hinder a wider EV uptake. Some of these barriers include the EVs'  limited autonomy, sparsely and unevenly dispersed Charging Stations (CSs), long charging sessions, and long waiting times at public CSs. 
In recent years, extensive research has been done to address the operational challenges entailed by these barriers~\citep{Kucukoglu2021}.

In this paper, we focus on long-distance trips, where an EV has to travel from an origin to a destination in the shortest amount of time.
This setting is relevant for private EV drivers, in cases where the shortest path between the origin and the destination has energy requirements exceeding the EV’s autonomy.
To reach her destination, the EV driver must choose a set of CSs to visit, and determine the amount of energy to charge at each of them.
This setting entails various challenges specific to the routing of EVs.
The first challenge concerns the limited EV autonomy.
In addition, despite global efforts supporting the increasing  adoption of EVs, the charging infrastructure is still sparse compared to the network of Conventional Vehicle (CV) refueling stations. 
Hence, EVs are often forced to detour from their shortest path to recharge.
Furthermore, in contrast to CVs, where refueling time is very short and can be considered constant, the EV charging time is generally much longer, depends on the CS technology, and follows a non-linear function. 
These factors significantly complicate decisions related to the amount of energy to charge at each CS visit.

Logistics operators using EVs generally prefer using private CSs~\citep{Morganti2018}. By the term private, we mean that CSs are located at the operator's facilities or at those of their customers. In such cases, logistics operators may plan charging activities in a way that avoids conflicts due to several EVs using the same charging infrastructure~\citep[e.g.,][]{Froger2022}.
Conversely, privately owned EVs, even if fully charged at their origin, are bound to public CSs when traveling long distances. In such contexts, CSs may be busy upon the EV arrival, thus entailing uncertain waiting times before the charging operation can begin. We focus on the associated challenges of using  public CSs for long EV trips.

The relatively long charging sessions of EVs at public CSs imply potentially long waiting times. 
Therefore, estimating and accounting for waiting times when planning routes is fundamental. 
Such estimations are compounded by the fact that the charging times are uncertain.
Therefore, even when the number of EVs in the queue at a CS is known, the resulting waiting time may vary significantly.
In the existing literature, for example~\cite{Sweda2017} and~\cite{Kullman2021}, it is generally assumed that the actual number of EVs at a CS is realized only when the EV physically observes the queue of the CS.
However, recent technologies provide real-time binary information regarding the status of public CSs, i.e., busy or not.
For example, the mobile application/web service ``ChargePoint"~\citep{ChargePoint} provides publicly available real-time occupancy information of CSs in the US, Canada, Europe, and Australia.
We refer to this binary information as an \textit{Occupancy Indicator} (OI). 
We note that, as this indicator is binary, it does not provide information regarding the number of EVs in the queue at the CS. 

The objective of this paper is to integrate OI information in route planning and quantify their added-value.
To this purpose, we consider the problem of routing an EV from a given origin to a given destination in the least possible time. The time to reach the destination includes driving, waiting, and charging at CSs. We consider that the EV continuously receives OI updates related to all CSs in a given network. At each such update, we determine the CSs to visit along with their associated charging quantities.
We name the resulting problem as the Electric Vehicle Shortest Path Problem with charging station Occupancy Indicator information (\evspp{}). 
The \evspp{} belongs to the class of problems known as Electric Vehicle Shortest Path Problems (EVSPPs). Furthermore, we consider a number of CS charging technologies, each with a non-linear charging function, and model each CS by an  appropriately defined queue. Thus, the \evspp{} can be seen as a variant of the Resource Constraint Shortest Path Problem (RCSPP) with several layers of additional complexity, including uncertain dimensions. 

We formalize the \evspp{} as a Markov Decision Process (MDP).
The MDP formulation accounts for dynamic updates of the state of the system, such as OI changes or queue observations, and enables us to model the desired EV routing and charging policy.
We solve the \evspp{} by a reoptimization algorithm. At each system update, the algorithm computes a sequence of CS visits and corresponding charge quantities. The computed sequence is followed until the next system update is received, and the process is then iterated. The proposed algorithm consists of several components.
During the execution of the trip, at each iteration, when a CS status changes or when the actual queue is observed, the algorithm estimates the expected waiting time at each of the CSs. 
To do so we address two main challenges. The first challenge is  due to the stochastic arrival and service times of exogenous EVs at CSs. To address this challenge, in a preprocessing phase, we devise a simulation-based procedure to estimate the waiting time at each CS as a function of the EV arrival time and the elapsed time from the last queue observation or OI information change.
Since the EV may need to visit, charge, and potentially wait at intermediate nodes before arriving at a CS, estimating the EV arrival time
at a given CS is challenging. To address this second challenge, we propose an efficient heuristic algorithm that estimates a lower bound on the EV arrival time at each CS. This lower bound is computed based on the queue-related characteristics of the CS. Once waiting times are estimated, the algorithm then computes the sequence of CS visits and recharging quantities by converting the \evspp{} into a deterministic EVSPP. This is solved to optimality by the procedure proposed in \cite{Kullman2021a}.

We evaluate the performance of the proposed algorithm by conducting a comprehensive computational study. 
We designed six experimental settings on a large set of artificial instances. 
We then benchmark the proposed algorithm against an alternative algorithm that makes no use of the OI.
We demonstrate that, on average, the proposed algorithm outperforms the benchmark in terms of the total trip duration by 1.4\%-18.5\%. 
This is primarily due to waiting times reductions of 23.7\%-95.4\%.
Moreover, we perform sensitivity analyses with respect to changes of several parameters, such as the density of CSs in the region, the technology of CSs, and others.

Summarizing, the main contributions of the paper are as follows:
\begin{itemize}
    \item We introduce a new dynamic and stochastic variant of the EVSPP with real-time information on the status of CSs.
    The proposed problem accounts for several challenging aspects, including the fact that the EV is allowed to partially charge its battery, potentially multiple times, and charging times are modeled via state-of-the-art piecewise linear charging functions that depend on the CS technology.
    \item We propose an MDP formulation for the \evspp{} in order to model the desired EV routing and charging policy while accounting for dynamic updates of the system.
    \item We propose a reoptimization algorithm that optimizes the sequence of CS visits and charging amounts based on system updates. 
    The algorithm includes several innovative aspects, including a simulation-based preprocessing phase to approximate the dynamics of each CS, and a waiting time estimation technique based on approximated lower bounds for the EV arrival time at a given CS. 
    The combination of these two elements allows us to take advantage of an existing state-of-the-art solver to efficiently update the CS visit sequence at each OI update.
    \item We build a set of artificial instances and conduct comprehensive computational experiments. We show that our algorithm clearly outperforms the benchmark algorithm, which does not take advantage of OI information. 
    Furthermore, we identify several contexts in which the edge of the proposed method with respect to the benchmark is stronger.
\end{itemize}

The remainder of the paper is organized as follows. In Section~\ref{sec:literature}, we present the related literature.
We describe the \evspp{} and its MDP model in Section~\ref{sec:problem}. We present the proposed algorithm in Section~\ref{sec:solution} and show the computational results in Section~\ref{sec:results}.
Finally, we provide our concluding remarks in Section~\ref{sec:conclusion}.

\section{Literature review}\label{sec:literature}
Over the past decade, extensive research has been conducted on routing problems arising in the context of operating EVs. For a comprehensive review, readers are referred to the surveys of~\citet{Erdelic2019} and \citet{Kucukoglu2021}.
Contributions in this area can be classified under two categories. 
The first consists of vehicle routing problems tailored to EVs (EVRPs). These generally aim at routing a fleet of EVs to serve a set of customers~\citep[e.g.,][]{Froger2019, Kullman2021}. 
The second category deals with EV shortest path problems (EVSPPs), where the objective is finding the shortest path between a given origin-destination pair~\citep[e.g.,][]{Sweda2017, Baum2020}. 
Due to the EV limited autonomy, in both groups, CS visits are typically planned en route.
The proposed \evspp{} falls under the second category as it deals with an EVSPP.

Various objective functions have been considered in the literature for EVSPPs, such as minimizing the trip duration~\citep{Baum2019}, the energy consumption~\citep{Baum2021}, or a combination of costs~\citep{Sweda2017}. 
The \evspp{} minimizes the total trip duration. 
Given that we consider long-distance trips, the EV may need to recharge its battery multiple times during its trip.
Depending on the problem setting, different EV charging policies have been considered in the literature.
For example, \cite{Erdogan2012} and \cite{Andelmin2017} consider the Green VRP, in which an EV is fully charged in constant time upon visiting a CS.
This assumption fits the case of alternative fuel vehicles or the case of EVs with battery swapping stations.
As battery swapping is not highly common in practice, charging times at CSs are often considered to be dependent on the amount of energy to charge. 
The relationship describing the time needed to charge a certain amount of energy is captured by charging functions. 
Several contributions have considered linear charging functions, e.g., \cite{Felipe2014} and \cite{Desaulniers2016}. 
However, in practice, charging functions are not linear \citep{Pelletier2017}. 
Accordingly, \cite{Montoya2017} and \cite{Froger2019} propose a piecewise linear function to approximate the realistic charging curves. 
Indeed, the adoption of piecewise linear functions to model EV charging time is common in practice.
The interested reader is referred to \cite{OpenEV2022}, where such functions are reported for about 300 EV models.
Notably, the charging functions are also heavily dependent upon the used charging technology (i.e., the speed of the chargers). 
Therefore, several authors have considered different CS technologies, each coupled with a distinct charging function \citep[e.g.,][]{Montoya2017, Felipe2014}. 
In the \evspp{}, we consider multiple charging technologies and piecewise linear charging functions. We also consider partial charging, i.e., the amount of energy to be charged at each CS is part of the decision.

Based on their ownership, CSs can be classified into three categories: private, semi-private, and public.
In problems with private CSs~\citep[see, for example][]{Froger2022}, the planner has the full control over CSs.
Semi-private CSs are either shared by a set of companies~\citep{Koc2019} or they are public but reserved to companies for specific time windows \citep{Bruglieri2019a}. 
Since in these two cases the planner has control over CSs, waiting times can be avoided or accounted for in deterministic models.
In contrast, public CSs are subject to exogenous demands, whose corresponding charging times are unknown to the route planner.  
In the context of EVSPPs, this leads to uncertain waiting times when visiting a CS.

Different approaches have been adopted in the literature to model the waiting time at public CSs.
The first one assumes that the waiting time is deterministic but time-dependent~\citep[see, for example][]{Gareau2019}.
The second one explicitly considers stochastic waiting times. \cite{Sweda2017} and~\cite{Keskin2021} assume that waiting times follow a given probability distribution. 
However, they assume that the actual waiting time is realized when the EV arrives at a CS.
This assumption is a simplification of reality, as when an EV arrives at a CS it only observes the length of the queue. 
Since the service time of each vehicle is stochastic and vehicles in a queue do not usually communicate, an estimation of the  waiting time is still needed.
Assuming that the length of the queue is observed upon arrival, \cite{Kullman2021} proposed a different strategy for estimating the waiting times.
They modeled the CS as a queue and leveraged the specific characteristics of queues. In particular, they stored previous observations of CS queue lengths and used them to simulate the waiting times and estimate their expected values. 
In our paper we adopt a similar technique to describe the transient behavior of queues at CSs following direct queues observations or OI information updates. 
Other approaches utilizing queuing characteristics of CSs are also studied in~\cite{Keskin2019} and~\cite{Setak2019}. 
These contributions used the steady-state conditions of a CS queue to estimate its associated waiting time.

In contrast to the previously discussed papers, where information about CSs is solely gathered through a physical visit, in this paper, we take into account an additional source of information. In particular, we assume the knowledge of dynamically-updated CS occupancy information. Indeed, due to a number of technological advancements, such information is now available to  EV users through mobile applications \citep[e.g.,][]{ChargePoint}. As previously mentioned, we formalize such information via the OI information.  We note that the current use of this information is limited to push notifications while navigating to a user-selected CS \citep[e.g.,][]{Evway}.
We examine how this information can be exploited to dynamically reoptimize the route as new information arrives during the trip.

Two streams of solution methods are relevant for our problem. 
The first stream consists of solution methods for EVSPPs.
For example, \cite{Sweda2017} formulated an EVSPP with non-linear charging functions, where partial charging is allowed, using dynamic programming.
Specifically, they examined two models: the first focuses on finding an optimal charging solution for a given fixed path (i.e., to decide how much to charge at each visit), and the second model considers finding an optimal solution for both routing and charging.
For an EVSPP with piecewise linear charging functions and partial charging, \cite{Baum2019} develop exact algorithms based on modified versions of A$^*$ and contraction hierarchies, as well as heuristic variants of the proposed algorithms. 
We note that, in contrast to \cite{Sweda2017}, \cite{Baum2019} do not take into consideration stochastic waiting time.
As a subcategory of EVSPPs, \cite{Wang2020} and \cite{Guillet2022} studied a problem of dynamically assigning an available CS to the EV upon receiving a charging operation request from the EV driver. The amount of charge is not considered as a decision factor in these problems.

The second stream consists of solution methods developed for Fixed Route Vehicle Charging Problems (FRVCPs). 
We refer interested readers to~\cite{Froger2019} for a complete review.
The FRVCP, introduced by~\cite{Montoya2016}, is a sub-problem of the EVRP. 
Given a fixed sequence of customers, the FRVCP determines the sequence of CSs to visit and establishes the amount of energy an EV should charge at each of them. 
The objective of the FRVCP is to minimize the completion time of the route, while ensuring the route feasibility with respect to the battery capacity.
Considering a linear charging function and partial charging, \cite{Roberti2016} proposed a labeling algorithm for the FRVCP.
Considering a piecewise linear charging function, different CS technologies, and partial charging, \cite{Montoya2017} proposed a mixed integer linear programming model, as well as a hybrid metaheuristic combining Iterated Local Search with Heuristic Concentration. 
\cite{Froger2019} proposed an exact labeling algorithm for the same problem considered in \cite{Montoya2017}.
We note that in this second stream, all problems consider private CSs and therefore do not account for stochastic waiting times.
In this study, once we estimated the waiting times, we adopt the labeling algorithm proposed by~\cite{Froger2019} to solve the deterministic version of our problem.

In summary, this paper examines the EVSPP-OI, which is a complex variant of the EVSPP considering public CSs, multiple charging technologies, piecewise linear charging functions, and a partial charging policy. In addition, we consider that the EV generally experiences stochastic waiting time at CSs.
Furthermore, the access to OI information adds a dynamic component to the problem, as routes can be modified in real-time based on newly received information.
To the best of our knowledge, the EVSPP-OI has not been addressed in the literature, and existing methods cannot be easily adapted to handle it.
This constitutes the main motivation of the present paper.

\section{Problem definition}\label{sec:problem}
In Section~\ref{sec:problemdescription}, we formally define the problem and formulate it in terms of an MDP in Section~\ref{sec:mdp}.

\subsection{Problem description}\label{sec:problemdescription}
We consider an EV that must travel from an origin $o$ to a destination $d$. 
A set of $\mathcal{C}=\{1, ..., n\}$ of public CSs may be used by the EV to charge its battery. 
Each CS $i\in \mathcal{C}$ has one charger of type $f_i\in F$, where $F$ is the set of available charging technologies. Specifically, we consider slow, fast, and normal chargers.

Since the considered CSs are public, we account for the presence of other vehicles, seeking to charge their batteries.
In the rest of the paper, we refer to the specific EV for which we plan the trip as \textit{the EV}, while the other vehicles are referred to as \textit{exogenous EVs}.
Similar to~\cite{Kullman2021}, we model each CS $i \in\mathcal{C}$ as an $M/M/1/\kappa$ queue, where the first two components refer to Poisson distributed arrival rates and service rates (i.e., exponential inter-arrival times and service times), and the third component indicates that the queue has one server. Finally, the last component $\kappa$ is the capacity of the queue at a CS. 
For example, $\kappa=1$ entails that exogenous EVs cannot wait at the CS, whereas $\kappa=2$ implies that at most one exogenous EV can wait at a CS while another exogenous EV is being charged. 
Thus, in general, exogenous EVs do not wait at a busy CS if $\kappa - 1$ EVs are already waiting in the queue. 
Restricting the capacity of the queue stems from the fact that a limited number of parking spots are typically available at each CS for waiting EVs.
The CS-related queuing parameters are defined as follows:  the arrival rate of the exogenous EVs at CS $i\in \mathcal{C}$ follows a Poisson distribution with expected value $\lambda_i$ {(EVs/hour)}, and the service rate (EVs/hour) follows a Poisson distribution with expected value $\mu_i$.  
The service rate of CS $i\in \mathcal{C}$ exclusively depends on its technology $f_i\in F$.

\begin{figure}
    \centering
    \includegraphics[width=0.5\textwidth]{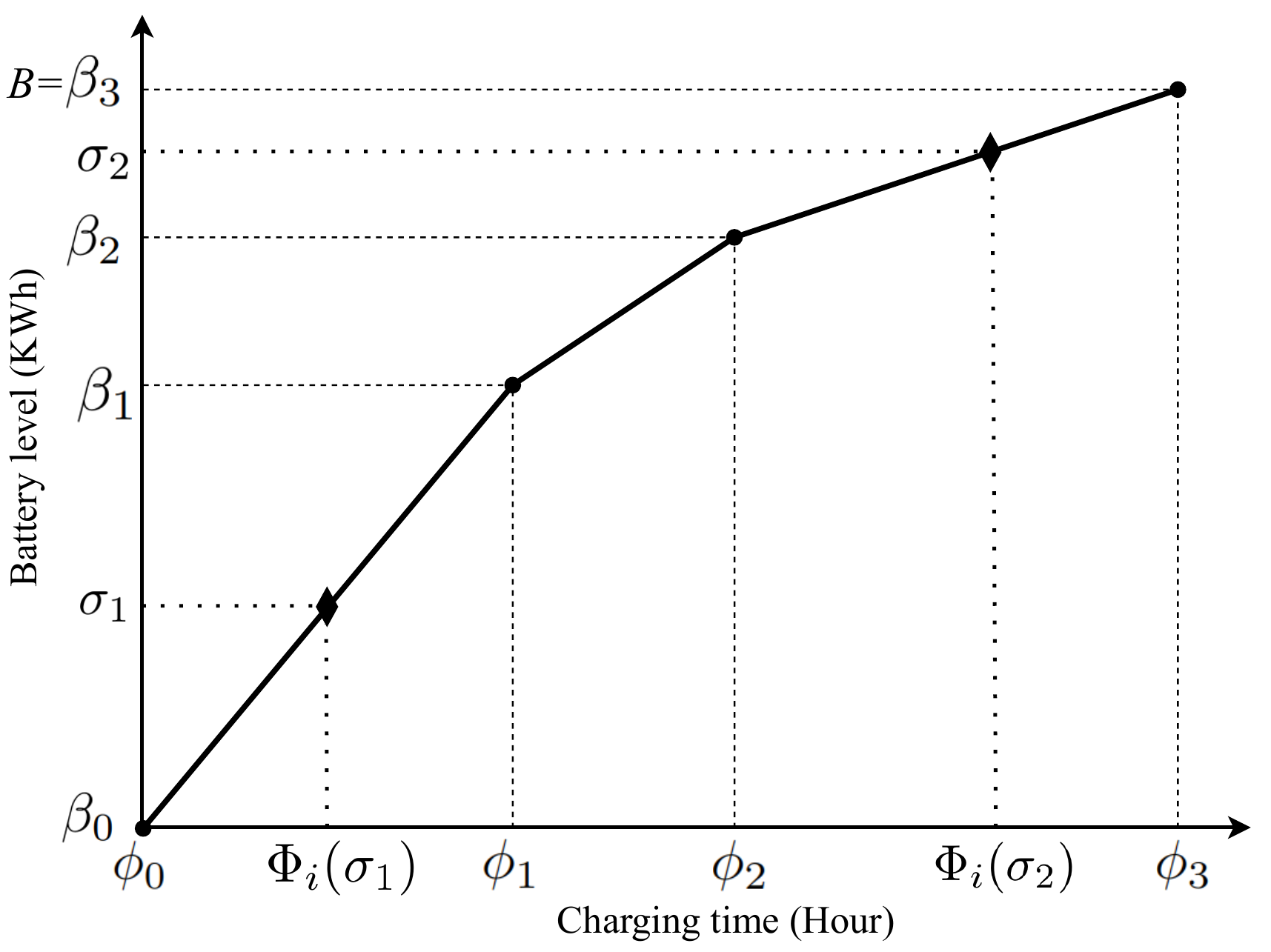}
    \caption{Piecewise linear approximation for charging function in CS $i$. Adapted from \citet{Froger2019}}
    \label{fig:curve}
\end{figure}

Let $G(N, A)$ be a directed graph, where $N=\mathcal{C}\cup \{o, d\}$ is the set of nodes and $A=\{(i,j)|i,j\in N, i\neq d, j\neq o\}$ is the set of arcs. Each arc represents the fastest path between a pair of nodes in the road network. 
In particular, we define $t_{ij}$ and $e_{ij}$ as the travel time and the energy consumption of the arc $(i,j)\in A$. Both the travel time and the energy consumption follow the triangular inequality. The EV has a battery with capacity $B$.
As suggested by~\citet{Montoya2017}, we approximate the charging time at each CS by a piecewise linear function.
In particular, each CS $i$ is associated with a piecewise linear concave charging function $\Phi_i(\sigma)$ that returns the charging time required to charge the EV from the battery level 0 to $\sigma$.
Let $z_i$ be the number of breakpoints of $\Phi_i(\sigma)$. 
Figure~\ref{fig:curve} illustrates an example of function $\Phi_i(\sigma)$ when $z_i=3$. 
In this figure, $\beta_z$ and $\phi_z$ represent the battery level and the charging time at breakpoint $z\in \lbrace 0,\dots,z_i\rbrace$, respectively.
We define function $C_i(\sigma_1, \sigma_2)=\Phi_i(\sigma_2) - \Phi_i(\sigma_1)$ to calculate the charging time required for the EV to charge at CS $i$ starting from battery level $\sigma_1$ until reaching battery level of $\sigma_2$, where $0\leq \sigma_1 \leq \sigma_2 \leq B$.

We focus on situations where $o$ and $d$ are relatively far from each other (i.e., $e_{od} \gg B$), so the EV requires recharging one or multiple times en route.
The EV is allowed to partially charge in the CSs. 
Consequently, the amount of charge at each charging session is a decision variable along with the sequence of CS visits.
We consider a random binary indicator for each CS $i\in\mathcal{C}$, denoted by $q_i(t)$, representing the real-time occupancy indicator of that CS at operating time $t$.
This indicator provides information about the occupancy of each CS. 
Specifically, when there are no EVs in the queue nor an EV that is being charged, then $q_i(t)=0$, otherwise $q_i(t)=1$.
In the latter case, the actual number of EVs in a CS is only realized when the EV physically arrives at the CS and observes the queue.
If the EV arrives at a busy CS, it may wait in the queue or leave for another CS.
As new information regarding the occupancy indicator becomes available at any time $t$, the estimated waiting times at all CSs will be updated and the EV may reoptimize the next node to visit and the corresponding amount of charge accordingly (we explain these procedures in Section~\ref{sec:solution}). 
However, we assume that once the EV starts charging, it will continue charging as planned irrespectively of new information being revealed during charging.
The \evspp{} aims at finding the optimal EV routing and charging policy in order to minimize the total duration of the trip from $o$ to $d$, which is comprised of waiting, driving, and charging times.

\subsection{MDP formulation}\label{sec:mdp}
Solving the \evspp{} requires to establish a policy determining the next node to visit, along with the target battery level to reach when leaving a CS node.
To this purpose, we formulate the \evspp{} as an MDP. We define a decision epoch as a point in time at which one of the following events occurs: (i) the EV leaves the origin, (ii) at least one occupancy indicator changes value, (iii) the EV arrives at a CS, (iv) the EV finishes waiting and the CS becomes available to serve it, or (v) the target battery level is reached and the EV leaves the CS.

We denote the state of the system at epoch $k$ as $s_k$. 
It is made up of three components: $(s^{v}_k, s^{c}_k, t_k)$, where $s^v_k$ is the state of the EV, $s^c_k$ is the state of the CSs, and $t_k$ is the time of the system at decision epoch $k$. 
More in detail, $s^{v}_k=(b_k, l_k, p_k)$ indicates the current battery level of the EV, its location, and its status, respectively. 
The status of the EV is either \textit{Driving}, \textit{Waiting}, or \textit{Charging}.
We represent the state of CSs as $s^{c}_k = [(q_{ik}, g_{ik})]_{i\in \mathcal{C}}$, where $q_{ik} = q_i(t_k)$ indicates if CS $i$ is occupied at epoch $k$ or not, and $g_{ik}$ is the elapsed time since the last time $q_{ik}$ has changed.
The terminal decision epoch $K$ is the point in time at which the EV arrives at the destination $d$.
Recalling that the graph $G(N, A)$ refers to the static position of CSs, origin, and destination, we define the graph $G_k(N_k, A_k)$ by including the node $c_k$, which indicates the current location of the EV at decision epoch $k$. 
In particular, when the EV is waiting or charging at a CS, then $c_k \in \mathcal{C}$, and $G_k=G$. 
However, when the EV is driving, $c_k$ is a new node that is added to $N_k$, i.e., $N_k=N\cup \{c_k\}$. In this case, we assume that all distances between $c$ and any node in $N$ can be retrieved, and thus $A_k=A\cup\{(c_k,i)|i\in N\setminus \{o\}\}$. For notational convenience, in what follows, we omit the index $k$ from $c_k$.

\begin{figure}[!htbp]
    \centering
    \includegraphics[width=.6\textwidth]{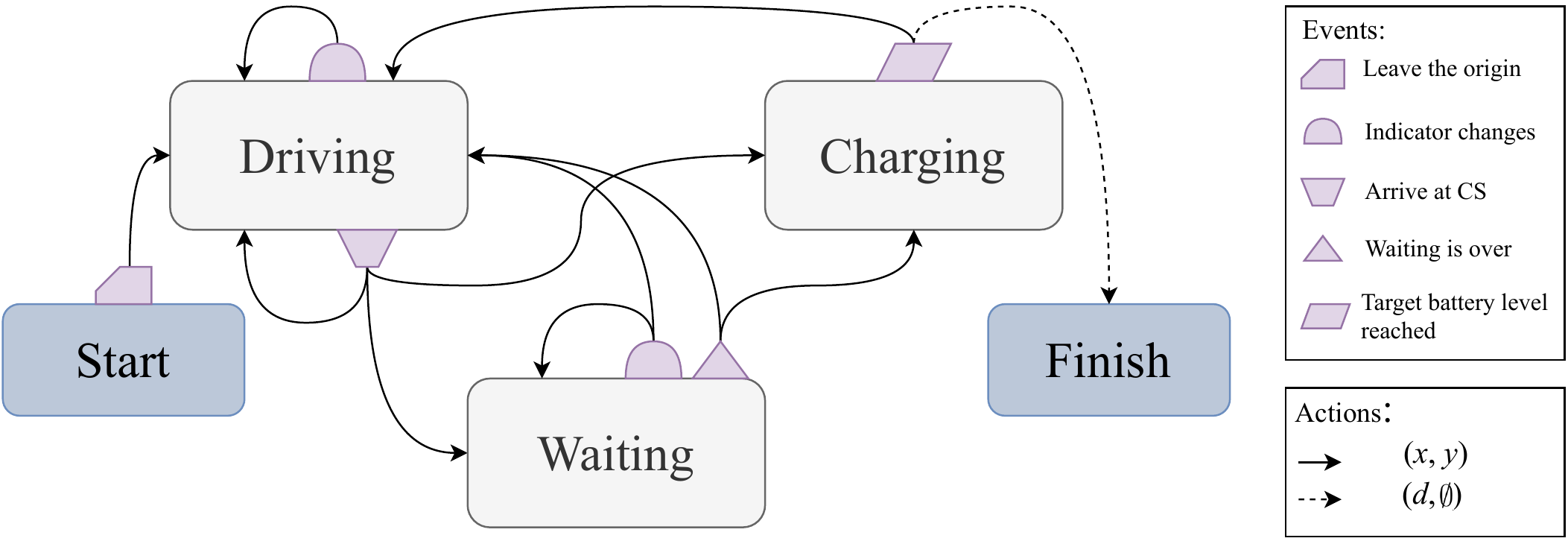}
    \caption{Dynamics of the MDP}
    \label{fig:mdp}
\end{figure}

At each decision epoch $k$, we define the action as a tuple $(x_k, y_k)$, where $x_k\in \mathcal{C}\cup \{d\}$ is the next node to visit and $y_k$ is the target battery level when leaving $x_k$. The action $(x_k, y_k)$ belongs to the action space $X_k$, defined as follows:
\begin{equation}
    X_k = \begin{cases}
    \bigl\{(d, \emptyset)\bigr\}, & e_{cd} \leq b_k \\
    \Bigl\{(x, y)|x \in \mathcal{C}, ~e_{cx} \leq b_k,~ y\in [b_k-e_{cx}, B]\Bigr\}, &  \textrm{otherwise}.
    \end{cases}
\end{equation}
Figure \ref{fig:mdp} illustrates the dynamics of the proposed MDP in the form of a Markov Chain. 
In this figure, arrows indicate actions, and small shapes (violet-colored) refer to the events that trigger the decision epoch. 
For example, when the EV is driving (i.e., $p_k=$ Driving) and arrives at a CS (i.e., event \includegraphics[width=15pt]{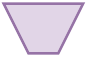} is triggered), it will start charging (i.e., $p_{k+1}=$ Charging) if the CS is available. Otherwise, the EV may choose to wait at the CS (i.e., $p_{k+1}=$ Waiting) or choose to travel to another CS (i.e., $p_{k+1}=$ Driving).
We note that while the EV is charging, ``Target battery level reached" is the only event that triggers a decision epoch. Therefore, no decision is made during the charging session.However, information regarding possible changes in the OI is always stored and accessible to the EV. 

The objective function of our problem is to minimize $t_K$. 
Accordingly, we define the cost of taking action $(x_k, y_k)$ in state $s_k$, denoted by $T(s_k, (x_k, y_k))$, as the time the EV spends between two decision epochs. Thus, we define the cost function as follows:
\begin{equation}
    T(s_k, (x_k, y_k))= t_{k+1} - t_k, ~\forall k\in \{0, 1, ..., K-1\}.
\end{equation}

\noindent Finally, we define the value of being in state $s_k$ using the Bellman equation as follows:
\begin{equation}
    V(s_k) = \min_{(x_k, y_k)\in X_k} \mathbb{E}_{\omega} \Bigl[T(s_k, (x_k, y_k)) + \gamma V(s_{k+1})\Bigr], ~\forall k\in \{0, 1, ..., K-1\},
\end{equation}
where $\omega$ refers to the exogenous information describing the stochastic arrival and service times at the CSs, and $\gamma<1$ is a discount factor to signify the importance of actions in the near future over farther ones.
According to the proposed formulation, we are looking for a policy that minimizes $V(s_k)$ for each  $s_k$.

\section{Solution method}\label{sec:solution}
The \evspp{} is challenging in several ways.
If we consider the simplified case where no recharging is allowed, the problem reduces to an RCSPP~\citep{Schiffer2018}, which is known to be $\mathcal{NP}$-hard~\citep{Garey1990}.
Including the recharging decisions in the \evspp{} significantly enlarges the decision space, when compared to a regular RCSPP. 
Moreover, the uncertainty of the arrival times and service times of the exogenous EVs at the CSs, which lead to stochastic waiting times, adds another layer of complexity to the problem.  

\begin{figure}[!htbp]
    \centering
    \includegraphics[width=0.8\linewidth]{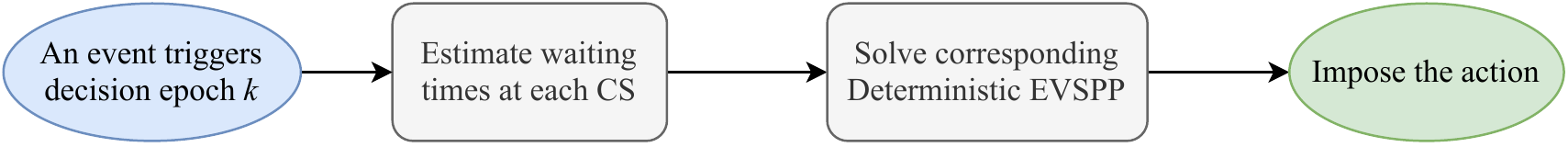}
    \caption{A schematic overview of \roi{}}
    \label{fig:reopt-overview}
\end{figure}

In order to solve the \evspp{}, we propose an online method that, at each decision, reoptimizes the next CS to be visited and the corresponding target battery level. As such, these decisions are reoptimized based on the new revealed information pertaining to the stochastic components of the problem. We call our method the Reoptimization method with Occupancy Indicator information (\roi).
Figure~\ref{fig:reopt-overview} shows an overview of \roi{} for each decision epoch. 
Accordingly, as an event triggers a decision epoch, we estimate the waiting time of the EV at each CS.
We then solve the corresponding deterministic EVSPP, where waiting times values are set to the computed estimates.

The first block in Figure~\ref{fig:reopt-overview} estimates the waiting times at all CSs.
The majority of studies in the literature modeling queues focus on estimating the expected waiting time under steady-state conditions. 
The steady-state of the queue is a unique state to which the queue converges after a long time, regardless of its observed state. 
Therefore, if the EV arrives a sufficiently long time after the last time the state of the queue was observed, the steady-state expected waiting time is a fairly accurate approximation of the actual waiting time.
However, in a relatively short period following a queue observation, the queue can be far from steady-state conditions, which makes the steady-state expected waiting time an inaccurate estimation of the actual waiting time.
In the \evspp{}, the OI provides real-time partial information regarding the occupancy of CSs, and thus the observed states of the CSs are frequently updated.
In such cases, the steady-state expected waiting time may be inaccurate, especially if the EV arrives at a CS shortly after the last time the state of the CS was observed.
We account for the information stemming from the OI in estimating the waiting times. 
This yields estimates that are closer to the actual waiting time than the steady-state expectation. 
Therefore, one of the main contributions of our work lies in introducing a waiting time estimation procedure that is compatible with the \evspp{} setting.

To determine the waiting time at each CS in a given decision epoch $k$ based on OI information, two main challenges need to be addressed. The first challenge arises due to the stochastic arrival and service times of exogenous EVs at CSs, making the waiting time subject to uncertainty. To address this challenge, we propose a method to estimate the expected waiting time at each CS $j$ using the OI information. As we model each CS as a queue of the form $M/M/1/\kappa$, we propose a simulation-based approach to estimate the waiting time at each CS. The resulting estimate depends on the arrival time of the EV at a CS. This leads us to the second challenge.
In particular, since the EV may need to visit, charge, and potentially wait at intermediate nodes in $\mathcal{C}\setminus j$ to travel from node $c$ to $j$, and given that the corresponding path (including the visit sequence and charging plan) is not known in advance, the arrival time of the EV at a given node $j$ may be challenging to compute. 
Even if we assume that the waiting time at each CS deterministically follows the proposed time-dependent function estimate, the resulting problem will be a time-dependent shortest path problem with constrained resources, some of which (specifically the charge of the battery) can be renewed according to multiple non-linear functions depending on the considered charging technologies. Developing an exact solution method for this extremely challenging problem is out of the scope of the present paper, especially considering the fact that we are aiming at developing a fast online algorithm.
We propose instead a fast heuristic that approximates a lower bound on the arrival time at each CS $j$ given the current location $c$. The proposed heuristic algorithm has a time complexity of $\mathcal{O}(m)$, where $m\leq|\{(i,j)|e_{ij}\leq B, i,j\in N_k\}|$.

Once waiting times are estimated, the \evspp{} is converted into a deterministic EVSPP (i.e., the second block in Figure \ref{fig:reopt-overview}). 
This problem entails deciding which CSs to visit and how much energy to charge at those CSs. 
To this purpose, we adopt the labeling algorithm proposed by \cite{Froger2019}, using the Python package of \citep{Kullman2021a}. The algorithm receives the battery state of the EV, travel time and energy consumption between nodes, and piecewise linear charging functions as input and returns a complete route and charging plan as the output, leading to the shortest total trip duration.
Since this algorithm does not specifically consider the waiting time in planning the route and the charging strategy, we use a modified version of the travel time, denoted by $\tau^k_{ij}$, as the input of the algorithm. As detailed in Section \ref{sec:waitingtime}, we define $\tau^k_{ij}$ as the sum of the direct travel time between two nodes $i$ and $j$ and the expected waiting time at CS $j$, given the state of the system at decision epoch $k$ and considering that the EV has passed $i$ immediately before $j$. 
Recalling that only an estimation of the waiting time at each CS can be obtained, we define $\Tilde{\tau}^k_{ij}$ as the approximation of the actual $\tau^k_{ij}$. We recursively apply this approximation on the path leading to $j$. 
Accordingly, while we keep the energy consumption values for traveling between nodes unmodified, we use $\Tilde{\tau}^k_{ij}$ as the travel time in the labeling algorithm, when reoptimizing the problem at decision epoch $k$.
For the sake of simplicity in notation, we will omit the index $k$ from $\tau^k_{ij}$ and $\Tilde{\tau}^k_{ij}$ in the following discussion.

In Section~\ref{sec:waitingtime}, we describe the procedure to estimate the waiting time used in computation of $\Tilde{\tau}_{ij}$. 
We outline the method used to solve the deterministic EVSPP in Section~\ref{sec:DSPP}.

\subsection{Waiting time estimation}\label{sec:waitingtime}
In this section, we develop a method to approximate the function that estimates the expected waiting time at a CS, assuming a given arrival time of the EV.
We denote this function by $W({q_{jk}}, {g_{jk}}, \Delta (\pi^c_j, t_k))$, which provides the estimated waiting time at CS $j\in \mathcal{C}$ computed at decision epoch $k$, given the states of all CSs, $(q_{jk}, g_{jk})$, and assuming that the EV arrives at $j$ after $\Delta (\pi^c_j, t_k)$ hours from $t_k$. 
The latter expression represents the time that the EV spends to travel from the current node $c$ to CS $j$ by taking an optimal path $\pi^c_j$ given the state of the system at time $t_k$. Thus, $\Delta (\pi^c_j, t_k)$ expresses the total driving, waiting, and charging times entailed by following an optimal path $\pi^c_j$. Formally, we define an optimal path $\pi^c_j=\{(c, b_k), (x^1,y^1), \dots, (x^{n'}, y^{n'}), (j,y)\}$ as a sequence of $(n'+2)$ CS-battery level pairs that the EV follows to travel from the current node $c$ to node $j$.
The first pair in this sequence refers to departing from the current node $c$ with the battery level of $b_k$.
Each pair of $(x^h, y^h)$, where $h\in\{1,\dots,n'\}$, indicates visiting CS $x^h\in\mathcal{C}$, possibly waiting for service, charging, and leaving that CS with battery level of $y^h\in[0, B]$. 
The last pair $(j,y)$ indicates arriving at node $j$ and planning to charge up to battery level of $y\in[0, B]$. 
If $e_{cj} \leq b_k$, we assume that no CS will be visited between $c$ and $j$ in $\pi^c_j$, and thus we set $\pi^c_j=\{(c, b_k), (j, y) \}$.
In general, $\pi^c_j$ must be energy feasible, i.e., $e_{c,x^1} \leq b_k, ~ e_{x^h, x^{h+1}} \leq y^h ~\forall h\in\{1,\dots, n'-1\}, {\text{ and }} e_{x^{n'}, j} \leq y^{n'}$.
For practical considerations, we assume that $\pi^c_j$ is elementary (i.e., a path that visits each node at most once).
For $\pi^c_j\ \neq \{(c, b_k), (j,y) \}$, we denote $i\in N$ the immediate predecessor of node $j$ in the path  (i.e., $i=x^{n'}$). 
As previously defined, let $\tau_{ij}({\pi^c_j})$ be the travel time from node $i$ to $j$ plus the expected waiting time at $j$ upon arrival, computed based on the state of the system at decision epoch $k$. Therefore,
\begin{equation}\label{eq:taugen}
     \tau_{ij}({\pi^c_j}) = t_{ij} + W\big({q_{jk}}, {g_{jk}}, \Delta (\pi^c_j, t_k)\big).
\end{equation}

As previously discussed, we face two challenges in estimating the waiting time function, $W({q_{jk}}, {g_{jk}}, \Delta (\pi^c_j,\allowbreak t_k))$, which we address by adopting two levels of approximation.
The first level of approximation, which is described in Section~\ref{sec:appfuncnn}, approximates the waiting time function $W(q_{jk}, g_{jk}, \Delta (\pi^c_j, t_k))$ for known input values.
The second involves the estimation of $\Delta (\pi^c_j, t_k)$. 
We note that while the optimal path $\pi^c_j$ is trivially computed for the cases where $e_{cj}\leq b_k$ (i.e., $\pi^c_j=\{(c, b_k), (j, y) \}$) and $e_{cd}\leq b_k$, computing an optimal path for all other cases is a complex task.
We propose a method to approximate $\Delta (\pi^c_j, t_k)$ using $\underline{\Delta}_{jk}$, which represents a lower bound on the duration of the trip between node $c$ and CS $j$ at epoch $k$. 
In Section~\ref{sec:estpath}, we describe the method providing an estimation of $\underline{\Delta}_{jk}$. 

\subsubsection{Approximating the function \texorpdfstring{$W(q_{jk}, g_{jk}, \Delta (\pi^c_j, t_k))$}{TEXT}}\label{sec:appfuncnn}
In this section, we propose a method to approximate the waiting time of the EV at CS $j\in \mathcal{C}$, for given input values $q_{jk}, g_{jk},$ and $\Delta (\pi^c_j, t_k)$.
In particular, for a CS $j$ with arrival rate $\lambda_j$, service rate $\mu_j$, and capacity $\kappa_j$, we approximate the expected waiting time in $\Delta (\pi^c_j, t_k)$ hours, knowing that $g_{jk}$ hours ago the state of the queue $q_{jk}$ changed from 0 to 1 (i.e., became occupied) or from 1 to 0 (i.e., became free).
Similar to \cite{Kullman2021}, we obtain the desired approximation by simulating the queue for a large set of scenarios, where each scenario represents the arrival and service time realizations of exogenous EVs. We observe the waiting times for each scenario. 
In particular, considering that the EV arrives at CS $j$ after $(g_{jk} + \Delta (\pi^c_j, t_k))$ hours have passed since the last change in $q_{jk}$, the corresponding queue needs to be simulated for $(g_{jk} + \Delta (\pi^c_j, t_k))$ hours, starting with $q_{jk}$ EVs in the queue. 

The function $W(q_{jk}, g_{jk}, \Delta (\pi^c_j, t_k))$ must be repeatedly calculated at each decision epoch in the reoptimization algorithm. 
However, running simulations online for a large set of arrival and service times every time the function is called is time-consuming.
One approach to reducing the load of online computations is to perform simulations in advance. Thus, we simulate a finite number of combinations of
$\lambda_j, \mu_j, \kappa_j$ corresponding to the considered queues in the system.
Specifically, as will be described in Section~\ref{sec:resultsinstances}, in our computational experiments, we consider 27 combinations of such parameters. 
Furthermore, by discretizing $(g_{jk} + \Delta (\pi^c_j, t_k))$ into steps of 0.1 hours, we execute the simulations for  each combination for 500K sets of exogenous arrival and service time realizations. 
In practical implementation, each simulation procedure entails simulating a queue of ($\lambda_j, \mu_j, \kappa_j$) for 20 hours, by starting from $q_{jk}$ EVs in the system, and recording the waiting time at discretized values of $(g_{jk} + \Delta (\pi^c_j, t_k))$.
We then store the results in a lookup table.
Retrieving waiting times from a lookup table, instead of executing simulations online, significantly reduces the computational time.

\begin{figure}[!htbp]
    \centering
    \includegraphics[width=.7\textwidth]{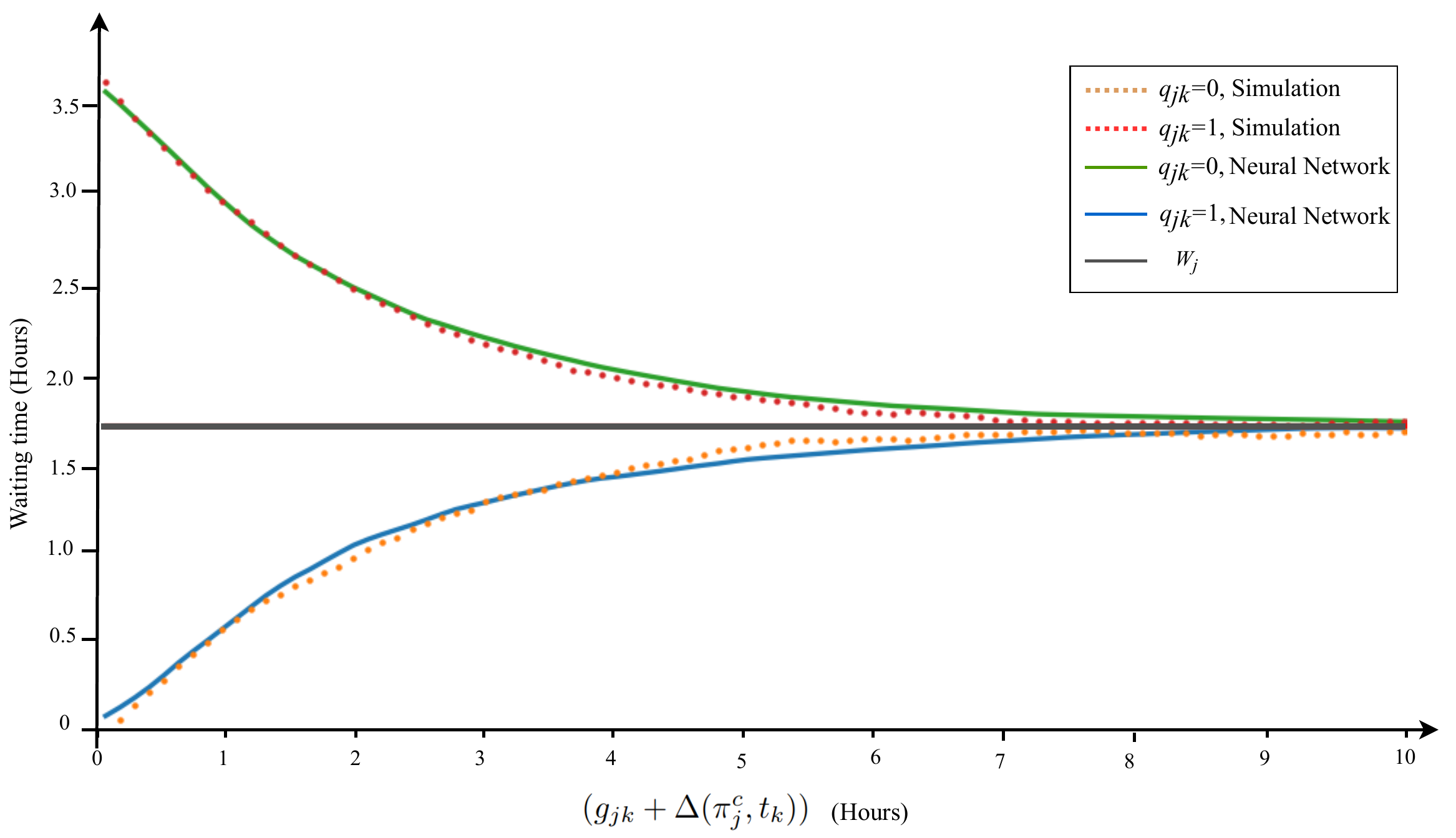}
    \caption{A demonstration for the performance of the trained approximate function}
    \label{fig:nnexample}
\end{figure}

The lookup table estimations can only be retrieved for discretized values of $(g_{jk} + \Delta (\pi^c_j, t_k))$. To further refine our estimations, we train a neural network that approximates the waiting time for continuous values of $(g_{jk} + \Delta (\pi^c_j, t_k))$. Specifically, we train a neural network on simulated waiting time series.
The proposed neural network has two hidden layers (with 256 and 128 neurons, each followed by a \textit{relu} activation function).
Figure~\ref{fig:nnexample} demonstrates the performance of the trained network for a queue at CS $j$ of the form $M/M/1/\kappa_j$, with $\lambda_j=0.25, \mu_j=0.28, \kappa_j=1$, compared with the simulated data.
In this figure, dotted curves (red and orange) show the trend of the estimated waiting times using simulation for a given arrival time. We note that each point in these two curves represents the estimated waiting time for an arrival time in $\{0, 0.1, \dots, 10.0\}$. 
Similarly, solid curves (green and blue) show the estimation obtained by the trained neural network. The horizontal (gray) line is the waiting time in the steady state, which we denote by $W_j$.

\subsubsection{Estimation of \texorpdfstring{$\underline{\Delta}_{jk}$}{TEXT}}\label{sec:estpath}
As described in Section~\ref{sec:waitingtime}, $\Delta (\pi^c_j, t_k)$ is equal to $t_{cj}$, when node $j$ is reachable from $c$ (i.e., $e_{cj}\leq b_k$).
In this section, we propose a method to estimate a lower bound on $\Delta (\pi^c_j, t_k)$, denoted $\underline{\Delta}_{jk}$, for $j\in\mathcal{C}$ when $e_{cj} > b_k$.

\begin{figure}
    \centering
    \includegraphics[width=0.7\textwidth]{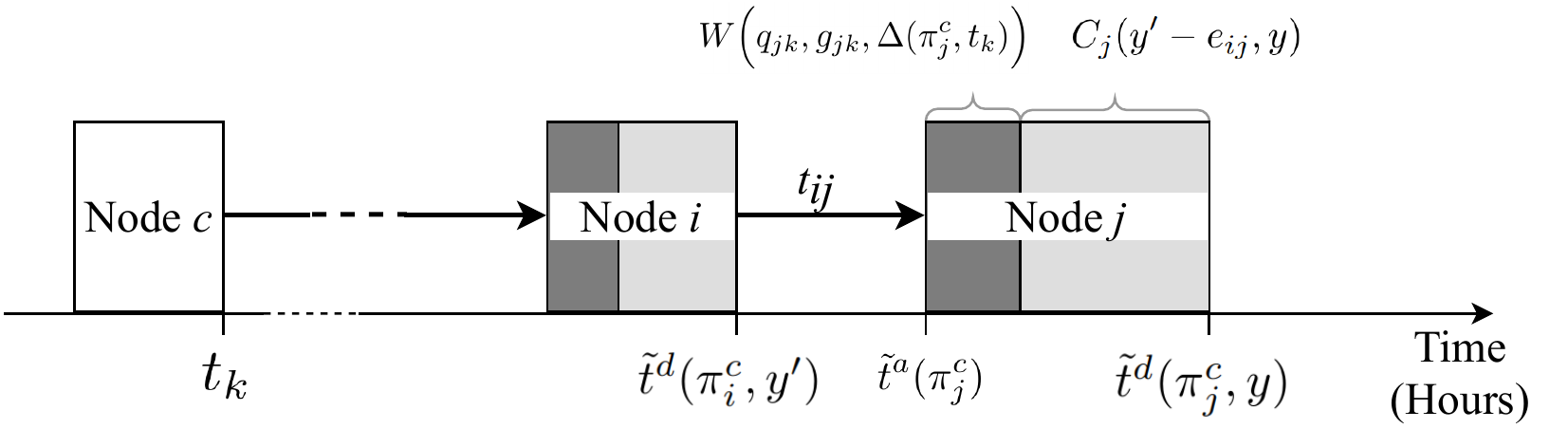}
    \caption{Demonstrating the recursive relationship between Equations~\eqref{eq:attime} and~\eqref{eq:TP0}}
    \label{fig:recursive}
\end{figure}

Let $\Tilde{t}^a(\pi^c_j)$ be the estimated arrival time of the EV at node $j$ by taking an optimal path $\pi^c_j$. Thus,
\begin{equation}\label{eq:deltaatime}
    \Delta (\pi^c_j, t_k)= \Tilde{t}^a(\pi^c_j) - t_k,
\end{equation}
Let $\Tilde{t}^d(\pi^c_j, y)$ be the estimated departure time of the EV from node $j$ after taking an optimal path $\pi^c_j$, including possible waiting at $j$, and charging up to the battery level $y$.
Given path $\pi^c_j=\{(c, b_k), (x^1,y^1), \dots, (x^{n'}, y^{n'}),\allowbreak (j,y)\}$, with the predecessor of node $j$ being $i\in N$ (i.e., $i=x^{n'}$), let the EV battery level when departing from $i$ be $y'$ (i.e., $y'=y^{n'}\in[e_{ij},B]$).
Accordingly, we define an optimal path $\pi^c_i$ such that $\pi^c_j=\pi^c_i\oplus (j, y)$. 
Therefore,
\begin{equation}\label{eq:attime}
    \Tilde{t}^a(\pi^c_j) = \Tilde{t}^d(\pi^c_i, y') + t_{ij}.
\end{equation}
We then express $\Tilde{t}^d(\pi^c_j, y)$ as follows:
\begin{equation}\label{eq:TP0}
    \Tilde{t}^d({\pi^c_j}, y) = \biggl[\Tilde{t}^a(\pi^c_j) + W\Bigl(q_{jk}, g_{jk}, \Delta(\pi^c_j, t_k)\Bigr) + C_j(y' - e_{ij}, y)\biggr].
\end{equation}
We assume that for $j=d$, $W(q_{jk}, g_{jk}, \Delta(\pi^c_j, t_k))=0$ and $C_j(.,.)=0$. 
Figure~\ref{fig:recursive} demonstrates the recursive relationship between Equations~\eqref{eq:attime} and~\eqref{eq:TP0}.

A lower bound on $\Delta(\pi^c_j, t_k)$ is computed by obtaining a lower bound on $\Tilde{t}^a(\pi^c_j)$, denoted $\Tilde{\underline{t}}^a_{jk}$.
Specifically, we define $\underline{\Delta}_{jk}=\Tilde{\underline{t}}^a_{jk} - t_k$.
Similarly, we let $\Tilde{\underline{t}}^d_{jk}$ be the lower bound on $\Tilde{t}^d (\pi^c_j, y)$. Then, by defining the set of nodes from which node $j$ is reachable as $N_k^j = \{i\in N\setminus\{j, d\}| e_{ij} \leq B \}$, we approximate Equations~\eqref{eq:attime} and~\eqref{eq:TP0} as follows:
\begin{equation}
    \Tilde{\underline{t}}^a_{jk} = \min_{i\in N_k^j} [\Tilde{\underline{t}}^d_{ik} + t_{ij}], ~\forall j \in N,
\end{equation}
\begin{equation}\label{eq:minapp}
    \Tilde{\underline{t}}^d_{jk}=\Bigl[\Tilde{\underline{t}}^a_{jk} + \underline{W}(q_{jk}, g_{jk}, \underline{\Delta}_{jk}) + \underline{C}_j\Bigr], ~\forall j\in N\setminus \{d\},
\end{equation}
where $\underline{C}_j$ and $\underline{W}(q_{jk}, g_{jk}, \underline{\Delta}_{jk})$ are the lower bounds on the charging time and the waiting time that the EV spends at CS $j$, respectively. 
Let $i^*_j=\arg\min_{i\in N_k^j} e_{ij}$, i.e.,  the node in $N_k^j$ from which traveling to node $j$ requires the least amount of energy.
We define $\underline{C}_j$ as follows:
\begin{equation}\label{eq:cmin}
    \underline{C}_j = C_j(0, \sigma^j_2 - \sigma^j_1),
\end{equation}
where $\sigma^j_1=B - e_{i^*_j,j}$ is the upper bound on the arrival battery level at $j$, and $\displaystyle{\sigma^j_2=\max(B - e_{i^*_j,j}, \min_{i\in N\setminus\{j, i^*_j\}} e_{ji})}$ is the lower bound on the target battery level when departing from $j$.  
This approximation assumes that the minimum amount of charging at CS $j$ happens when the EV arrives at CS $j$ with a maximum possible battery level (i.e., being fully charged at CS $i^*_j$) and charges enough to be able to travel to its closest node other than $i^*_j$. 
Given the shape of function $C_{j}(., .)$  (Figure~\ref{fig:curve}), it can be inferred that $C_j(\sigma_1, \sigma_2) \leq C_j(\sigma_1 + a, \sigma_2 + a)$, where $a\geq 0$. Therefore, in order to guarantee that $\underline{C}_j$ is a lower bound on the true value of $C_j(., .)$, we shift the charging range from ($\sigma^j_1$, $\sigma^j_2$) to $(0, \sigma^j_2 - \sigma^j_1)$. 

To compute $\underline{W}(q_{jk}, g_{jk}, \underline{\Delta}_{jk})$ as a lower bound on the actual waiting time upon the arrival of the EV at CS $j$, we use the characteristics of the approximate function $W(., ., .)$.
In general, once the approximate waiting time function for a given queue is developed (as detailed in Section~\ref{sec:appfuncnn}), its minimum point with respect to the arrival time, denoted $\underline{W}_j$, can be empirically obtained by inspection. We use $\underline{W}_j$ as an approximation for the lower bound of the waiting time upon arrival of the EV, i.e., $\underline{W}(q_{jk}, g_{jk}, \underline{\Delta}_{jk})$.

\begin{figure}[!htbp]
    \centering
    \includegraphics[width=0.5\textwidth]{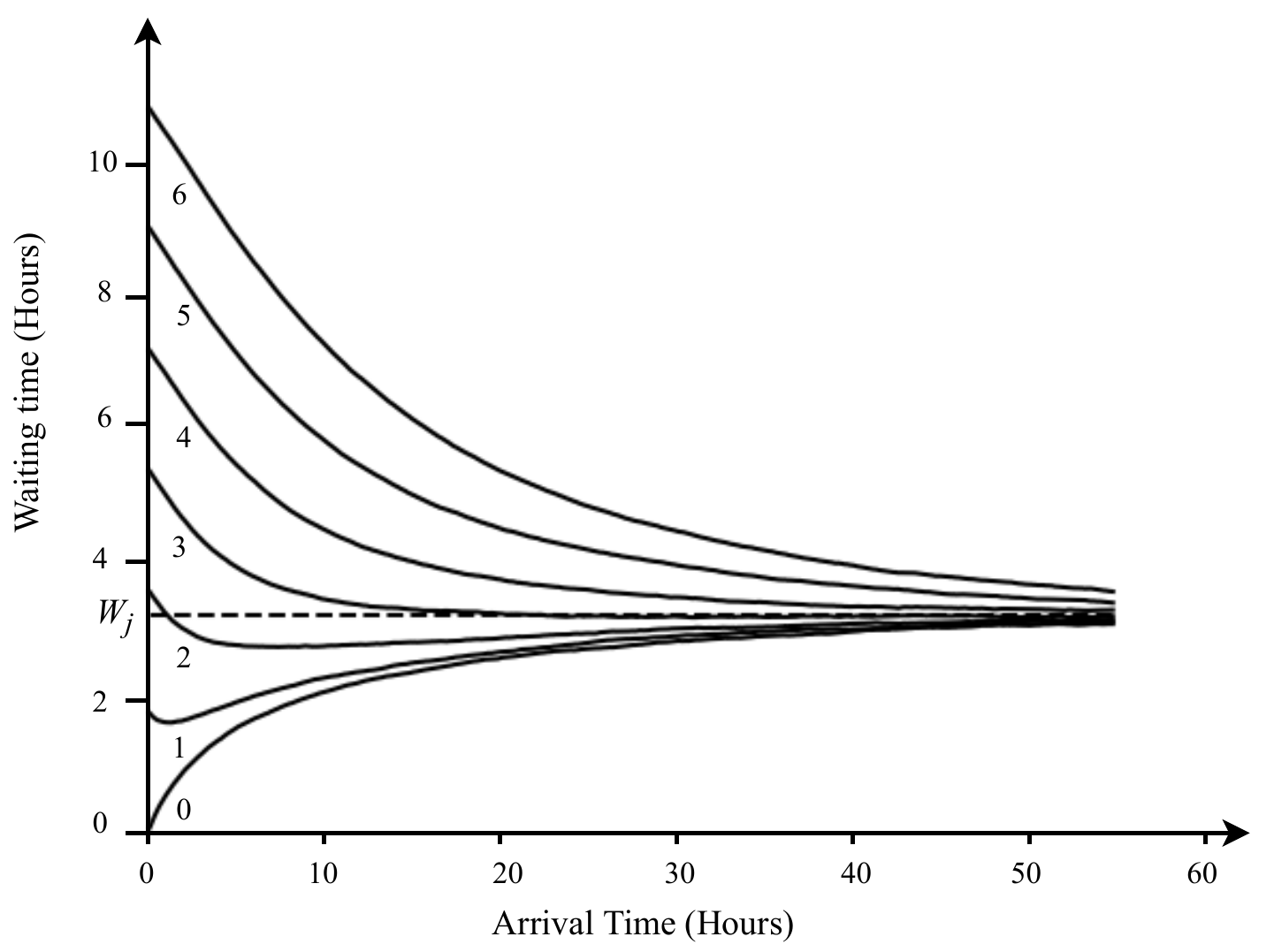}
    \caption{Steady-state convergence of the waiting time for a CS $j$ with a queue of the form $M/M/1/\infty$, where $\lambda_j=0.36$ and $\mu_j=0.56$ ($\rho_j=0.65$). \textit{Note}: The number written below each curve represents the queue length at time 0.}
    \label{fig:ss_gen}
\end{figure}

Based on our empirical analysis, we observed a systematic behavior of the waiting time function with respect to the arrival time. 
In some cases, this enabled us to infer a lower bound on $W(q_{jk}, g_{jk}, \underline{\Delta}_{jk})$ that is stronger than $\underline{W}_j$.
First, it is known that an $M/M/1/\kappa$ queue converges to a steady-state after a theoretically infinite amount of time, regardless of its initial state.
This behavior is illustrated in Figure~\ref{fig:ss_gen}. 
In this figure, we consider a queue with $\kappa=\infty$ to illustrate the steady-state convergence for various initial queue lengths (i.e., the number of EVs waiting for a CS plus the EV being served).
Within the context of our computational experiments, where $\kappa$ is considered to be less than or equal to three (described in Section~\ref{sec:resultsinstances}), we observed a comparable pattern of convergence, as depicted in Figure~\ref{fig:ss_gen}.
As can be seen from the figure, the average waiting time upon arrival converges to the steady-state waiting time $W_j$ regardless of the number of EVs at a CS at time 0.
The steady-state waiting time is defined as $W_j=\frac{L_j}{\mu_j}$, where $L_j$ is the queue length in its steady state. $L_j$ is computed as follows (interested readers are referred to~\cite{Thomopoulos2012} for details):
\begin{equation}\label{eq:lj}
    L_j = \rho_j(\frac{1-\rho_j^{\kappa_j}}{1-\rho_j^{\kappa_j+1}}) + \frac{\rho_j^2[1-\kappa_j\rho_j^{\kappa_j-1} + (\kappa_j-1)\rho_j^{\kappa_j}]}{(1-\rho_j)(1-\rho_j^{\kappa_j+1})},
\end{equation}
where $\rho_j = \lambda_j / \mu_j$ is the utilization rate of CS $j$. For example, for $\kappa_j=1$, $L_j=\frac{\rho_j}{1+\rho_j}$, and for $\kappa_j=2,~L_j=\frac{\rho_j(1 + 2\rho_j)(1-\rho_j)}{1-\rho_j^3}$.

\begin{figure}[!htbp]
    \centering
        \subfigure[Case A: $\lambda_j=0.73,~L_j=0.99,~W_j=0.89$]{
        \includegraphics[width=0.3\textwidth]{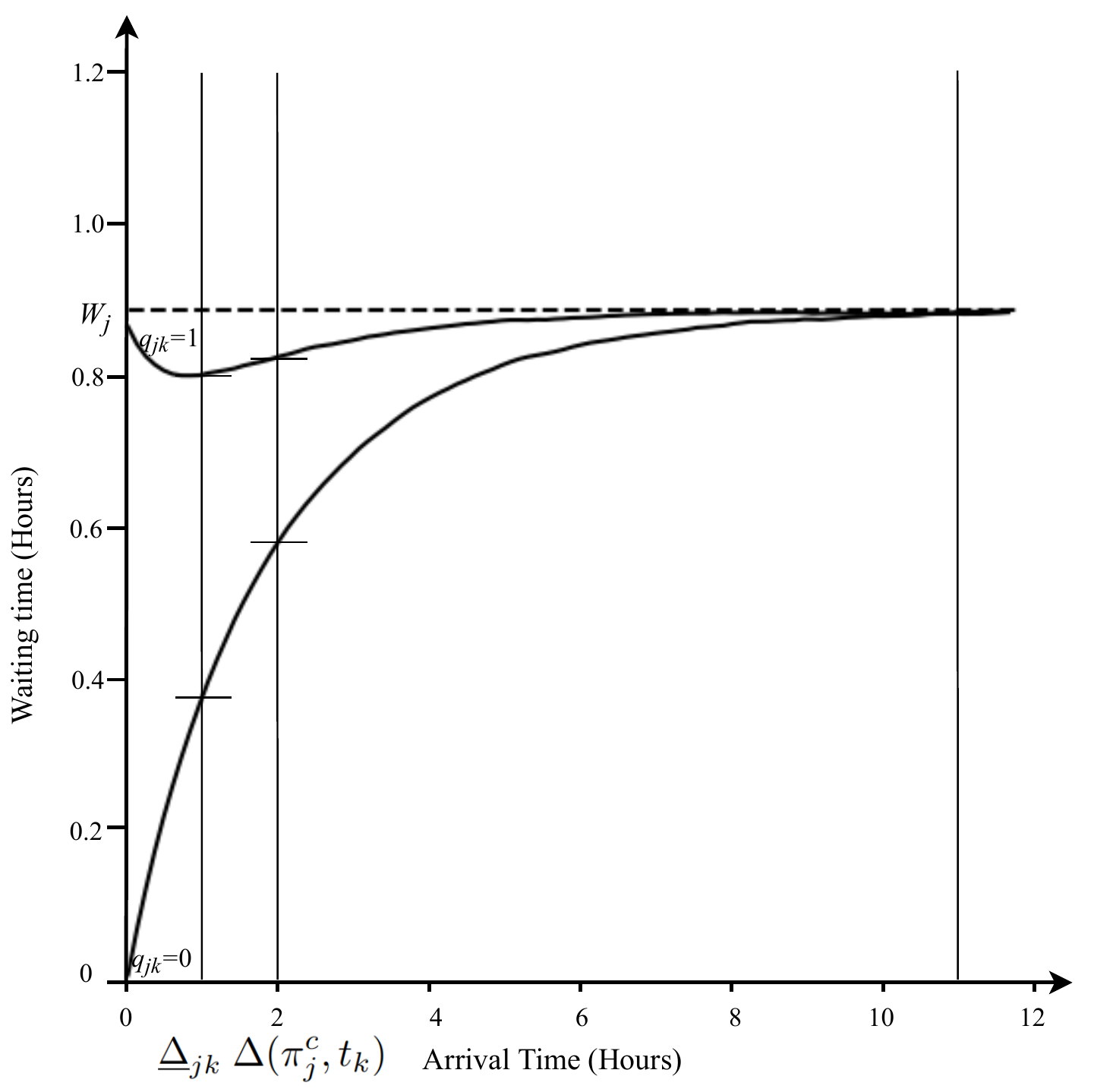}
            \label{fig:bounds_gen_A}
    }
    \subfigure[Case B: $\lambda_j=0.45,~L_j=0.57,~W_j=0.51$]{
    \includegraphics[width=0.3\textwidth]{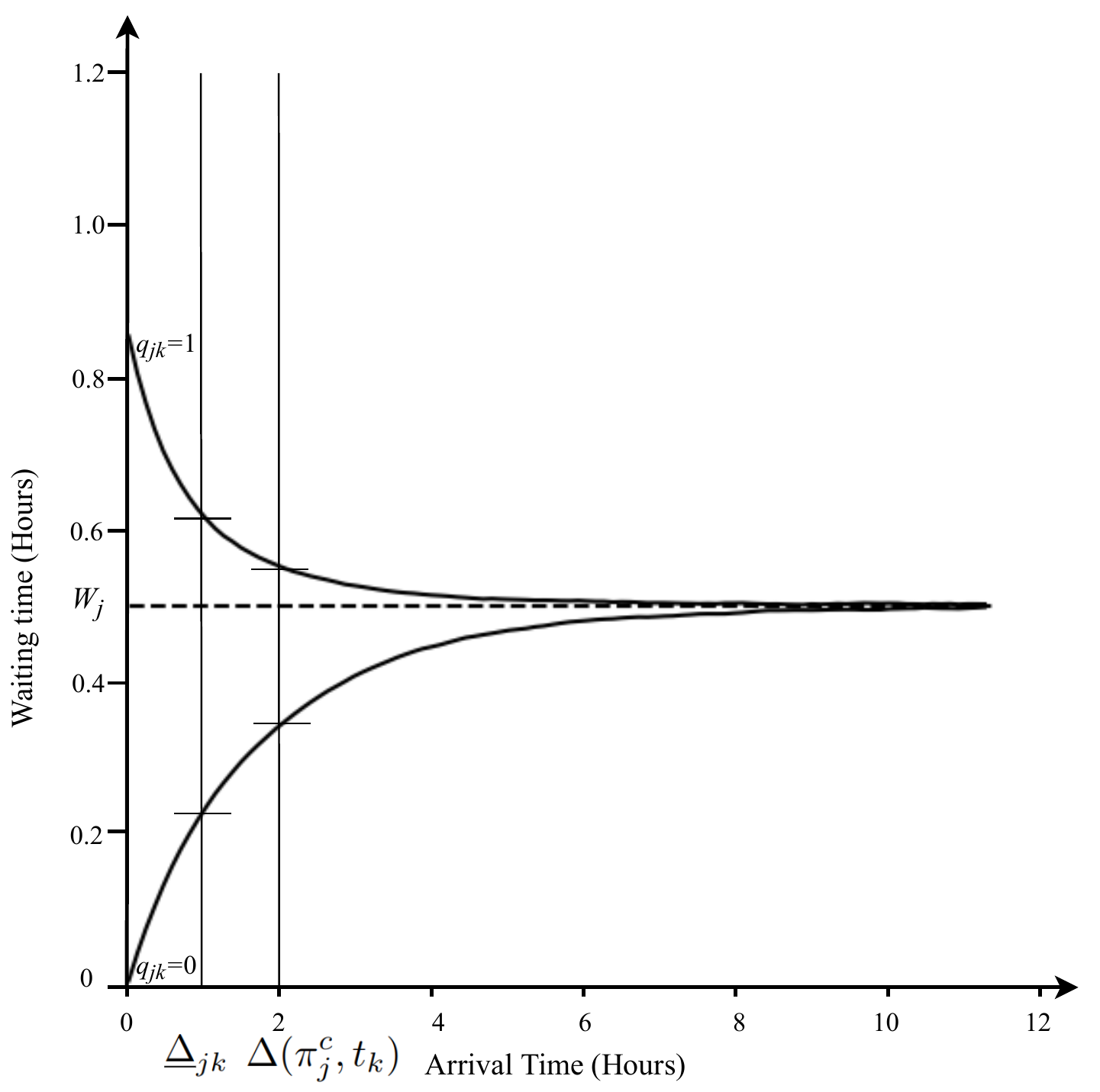}
        \label{fig:bounds_gen_B}
    }
    \subfigure[Case C: $\lambda_j=1.0,~L_j=1.36,~W_j=1.21$]{
    \includegraphics[width=0.3\textwidth]{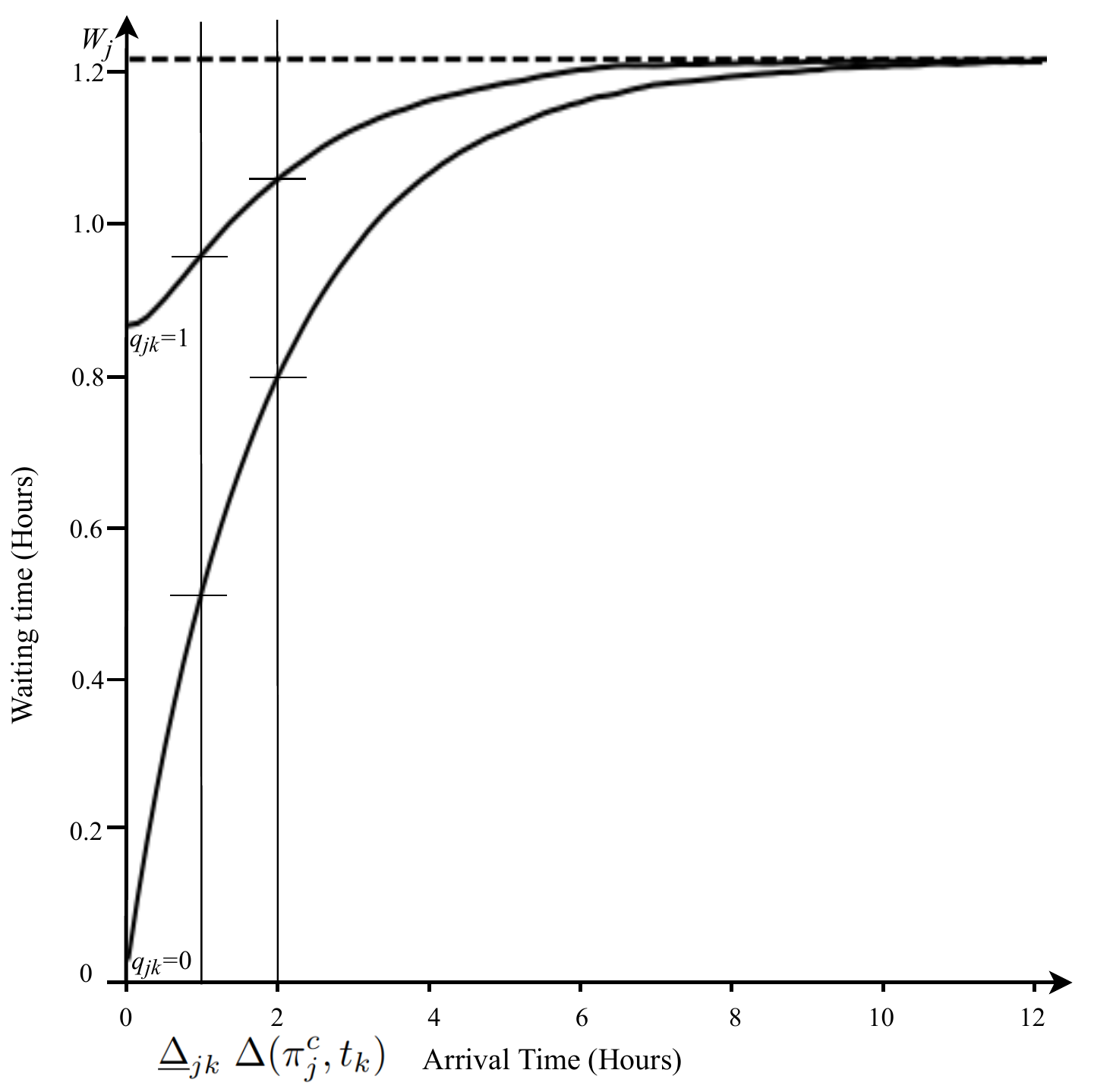}
        \label{fig:bounds_gen_C}
    }
    \caption{Approximate bounds for the waiting time for $\kappa_j=3$ and $\mu_j=1.12$}
    \label{fig:bounds_gen}
\end{figure}

We considered 27 combinations of $\kappa_i, \lambda_i$ and $\mu_i$ in our experiments (described in Section~\ref{sec:resultsinstances}). 
In those combinations, we observed three cases of waiting time behavior with respect to the arrival time.
Figure~\ref{fig:bounds_gen} illustrates one example of each case. 
Considering $\kappa_j=3$ and $\mu_j=1.12$, each sub-figure illustrates the simulated behavior of a given queue of the form $M/M/1/\kappa_j$ for two cases, $q_{jk}=0$ and $q_{jk}=1$, assuming that $q_{jk}$ has  just changed at epoch $k$ (i.e.,  $g_{jk} = 0$). We refer to Figures~\ref{fig:bounds_gen_A},~\ref{fig:bounds_gen_B}, and~\ref{fig:bounds_gen_C} as cases A, B, and C, respectively.   
The arrival rate $\lambda_j$ for the queues in these cases is 0.73, 0.45, and 1.0, respectively.
The curves indicated by $q_{jk}=1$ show that if the EV arrives at CS $j$ at time zero, the expected waiting time will be $\frac{1}{\mu_j}$ (0.89 in all  sub-figures). However, as the EV arrives later, the expected waiting time converges to $W_j$. 
As can be seen in all cases, for $q_{jk}=0$, we observe that $\underline{W}_j=0$. However, in these cases, the waiting time monotonically increases and converges to $W_j$. 
Therefore, when $q_{jk}=0$, $W(q_{jk}, g_{jk}, \underline{\Delta}_{jk})$ provides a stronger lower bound on the actual waiting time than $\underline{W}_j$. 
When $q_{jk}=1$, the waiting time function behaves differently for each case.
Specifically, in case A (with $L_j=0.99$), the waiting time is not monotonic with respect to the arrival time.
Therefore, we maintain $\underline{W}_j$ (equaling 0.8 in Figure~\ref{fig:bounds_gen_A}) as the lower bound on actual waiting time in this case.
In case B (with $L_j=0.57$), the waiting time function monotonically decreases and converges to $W_j$. Therefore, $W_j$ is our approximation for $\underline{W}(q_{jk}, g_{jk}, \underline{\Delta}_{jk})$ in this case. 
Finally, in case C (with $L_j=1.36$), the waiting time monotonically increases and converges to $W_j$. For this case, we set $\underline{W}(q_{jk}, g_{jk}, \underline{\Delta}_{jk})$ to $W(q_{jk}, g_{jk}, \underline{\Delta}_{jk})$. 
To summarize, we compute $\underline{W}(q_{jk}, g_{jk}, \underline{\Delta}_{jk})$ as follows:
\begin{equation} \label{eq:wt}
    \underline{W}(q_{jk}, g_{jk}, \underline{\Delta}_{jk}) = \begin{cases}
    W(q_{jk}, g_{jk}, \underline{\Delta}_{jk}) & q_{jk} = 0, \\
    \underline{W}_j & q_{jk} = 1\land \textrm{case A}, \\
    W_j & q_{jk} = 1\land \textrm{case B}, \\
    W(q_{jk}, g_{jk}, \underline{\Delta}_{jk}) & q_{jk} = 1 \land \textrm{case C}. \\
    \end{cases}
\end{equation}

In Algorithm 1, we present a forward recursive procedure  for computing $\Tilde{\underline{t}}^d_{jk}$ and subsequently $\Tilde{\underline{t}}^a_{jk}$ for all $j\in N$. 
The adjacency matrix $\mathcal{A}$ indicates if $j$ is directly reachable, when departing fully charged from $i$.
Specifically, we set $\mathcal{A}_{ij}=1$ if $(i=c \land e_{ij}\leq b_k)$ or $(i\neq c \land e_{ij}\leq B)$, otherwise we set $\mathcal{A}_{ij}=0$.
The algorithm starts by setting the earliest departure time from the current node $c$ to zero ($\underline{\Tilde{t}}^d_{ck} = 0$) and assigning $\underline{\Tilde{t}}^d_{jk} = \infty$ for all other nodes $j \in N$. At each step, it selects a node $i \in N_k$ with the minimum $\underline{\Tilde{t}}^d_{ik}$ that has not been selected before and updates $\underline{\Tilde{t}}^d_{jk}$ for all reachable nodes $j$ ($\mathcal{A}_{ij}=1$) from node $i$. The selected node $i$ is then added to a set $L$ to ensure it is not selected again.
This procedure terminates when $|L|=|N_k|$.

Once $\Tilde{\underline{t}}^d_{jk}$ is computed for all $j\in N$ using Algorithm~\ref{alg:cap}, we approximate $\tau_{ij} (\pi^c_j)$ by $\Tilde{\tau}_{ij}$, as follows:
\begin{equation}\label{eq:tauij}
    \Tilde{\tau}_{ij} = t_{ij} + W(q_{jk}, g_{jk}, \Tilde{\underline{t}}^d_{ik} + t_{ij} - t_k), {~\forall \{i,j\in N_k| e_{ij} \leq B \}}.
\end{equation}
By adopting Equation \eqref{eq:minapp} in Equation \eqref{eq:tauij}, we can express the relationship between $\Tilde{\tau}_{ij}$, the approximate function $W(.,.,.)$, and the estimated lower bound for arrival at CS $j$, as follows:
\begin{equation}\label{eq:tauij2}
    \Tilde{\tau}_{ij} = t_{ij} + W\Big(q_{jk}, g_{jk}, \big[\Tilde{\underline{t}}^a_{ik} + \underline{W}(q_{ik}, g_{ik}, \underline{\Delta}_{ik}) + \underline{C}_i + t_{ij} - t_k\big]\Big), {~\forall \{i,j\in N_k| e_{ij} \leq B \}}.
\end{equation}

We summarize the value of $\Tilde{\tau}_{ij}$ for all $i,j\in N_k$ as follows:
\begin{equation}\label{eq:tauijf}
    \Tilde{\tau}_{ij} = 
    \begin{cases}
    \infty & e_{ij} > B, \\
    \infty & i=j\neq c, \\
    \infty & i=j=c \land p_k\neq\textrm{Waiting}, \\
    W(q_{jk}, g_{jk}, 0) & i=j=c \land p_k=\textrm{Waiting} \land \kappa = 1, \\
    \xi_c/\mu_c & i=j=c \land p_k=\textrm{Waiting} \land \kappa > 1, \\
    t_{ij} + W(q_{jk}, g_{jk}, t_{ij}) & i=c \land j\neq c \land e_{ij} \leq b_k, \\
    t_{ij} + W(q_{jk}, g_{jk}, \Tilde{\underline{t}}^d_{ik} + t_{ij} - t_k) & \textrm{otherwise},
    \end{cases}, ~\forall i,j\in N_k.
\end{equation}
The first case in Equation~\eqref{eq:tauijf} sets $\Tilde{\tau}_{ij}=\infty$ if there is no energy feasible direct path from $i$ to $j$. 
The cases two to five indicate that the value of $\Tilde{\tau}_{ij}$ for two identical nodes (i.e., when $i=j$) is infinite, except in the situations where the EV is waiting at a CS (i.e., cases four and five). We note that these two cases occur either when the EV arrives at a CS or when it is waiting there.
In these cases, the EV observes the actual queue length as $\xi_j$.
When $\kappa > 1$, we compute the waiting time at $j$ as $\xi_j / \mu_j$. 
However, for $\kappa=1$, given that the OI provides perfect information regarding the beginning of charging operation for the exogenous EV at CS $j$, we use the waiting time estimation provided in Section~\ref{sec:appfuncnn}.
The sixth case relates to when directly traveling from $c$ to $j$ is energy feasible. To compute the waiting time for this case we consider $t_{ij}$ as the arrival time at $j$. 
In all other cases, we use the lower bounds computed for the arrival time at $j$ (i.e., $\Tilde{\underline{t}}^d_{ik}$) to estimate the waiting time used in $\Tilde{\tau}_{ij}$. 

\begin{algorithm}[!htbp]
\DontPrintSemicolon
Build the adjacency matrix $\mathcal{A}:\{0,1\}^{N_k\times N_k}$\;
Initialize $L=\{\}, ~\Tilde{\underline{t}}^d_{ck}=0, ~ \Tilde{\underline{t}}^d_{jk}=+\infty, ~\forall j\in N$\;

\While{$|L| \neq |N_k|$}
{
Choose $i$ such that $\Tilde{\underline{t}}^d_{ik} = \min \{\Tilde{\underline{t}}^d_{hk}|h \in N_k\setminus L\}$\;
$L:=L\cup i$\;
$L_i = \{j\in N| \mathcal{A}_{ij}=1 \}$\tcp*{Set of adjacent nodes to $i$}

\For{$j \in L_i$}{
$\underline{w} = \underline{W}(q_{jk}, g_{jk}, \Tilde{\underline{t}}^d_{ik} + t_{ij} - t_k)$\;
$\Tilde{\underline{t}}^d_{jk} := \min(\Tilde{\underline{t}}^d_{jk},~ \Tilde{\underline{t}}^d_{ik} + t_{ij} + \underline{w} + \underline{C}_j)$\;
}

}
\Return{$\Tilde{\underline{t}}^d_{jk},~ j\in N_k$}
\caption{Approximate $\Tilde{\underline{t}}^d_{jk},~ \forall j\in N_k$ at decision epoch $k$}\label{alg:cap}
\end{algorithm}

\subsection{Deterministic EVSPP} \label{sec:DSPP}
In Section~\ref{sec:waitingtime}, we described how to incorporate the OI information to estimate the waiting time in a CS upon the arrival of the EV. 
We  also defined $\tau_{ij}(\pi^c_j)$ as the direct travel time between $i$ and $j$ plus the waiting time at $j$ upon arrival.
We then estimated $\tau_{ij} (\pi^c_j)$ as $\Tilde{\tau}_{ij}$. 
At each epoch, we solve a deterministic EVSPP, where  we use the  estimated values of $\Tilde{\tau}_{ij}$ at each decision epoch as the deterministic travel time between two nodes, while keeping the energy consumption $e_{ij}$ between nodes unmodified. We use the output to establish the next CS to visit and the target battery level to charge at it.

To solve the deterministic EVSPP, we adopt the labeling algorithm proposed by \cite{Froger2019}. 
Given a fixed-route indicating the sequence of non-CS nodes to visit, with the objective of minimizing route completion time, this algorithm finds an optimal charging plan for the EV, while ensuring that route is energy feasible. 
In particular, given a set of CSs, this algorithm establishes which CSs must be visited between every two non-CS nodes in the fixed-route and decides how much to charge at each CS.
In our problem setting, this fixed-route is comprised of two nodes, $o$ and $d$ (i.e., the origin and the destination). For a detailed explanation, we refer the interested reader to \cite{Froger2019}. 
Briefly, the labeling algorithm assigns labels to the nodes in the graph $G_k$ tracking the time and charging level at each node. 
Piecewise linear charging functions are modeled similar to the ones used in our paper. 
In particular, each break-point of the piecewise linear charging function represents a possible charging interval. 
Therefore, for each break-point, a label indicating the departure time from that CS is generated. 
We note that this algorithm does not account for waiting times at CSs. 
To account for waiting times in the algorithm, we consider the estimated values of $\Tilde{\tau}_{ij}$ as the travel times between two nodes and maintain the energy consumption matrix $e_{ij}$ unchanged.

In order to make the labeling algorithm easier to use, \cite{Kullman2021a} implemented it as a Python package, the so-called \texttt{frvcpy}. 
This package receives the features of CSs (including their charging function, travel times between nodes, and energy consumption to travel between nodes) as input, along with the capacity of the EV and a fixed route, and returns an optimal energy-feasible route from $o$ to $d$, as a sequence of $n'+2$ CS-battery level pairs. Formally, it returns a path in the form of $\{(c, b_k),(x^1, y^1),\dots,(x^{n'}, y^{n'}), (d,\emptyset)\}$, where $x^h\in\mathcal{C}$ and $y^h\in[0,B], \forall h\in[1,n']$. 
We use this package in our {\roi} to solve the deterministic EVSPP.

\section{Computational results}\label{sec:results}
In this section, we present a series of computational experiments aimed at validating our solution method (i.e., \roi), and quantifying  the added-value of accounting for information provided by OI.
In Section \ref{sec:resultsinstances}, we describe the procedure to generate problem instances. 
In order to examine the added value of using the OI information, we benchmark the performance of our methodology against the setting where no such information is available. We describe this benchmark setting in Section~\ref{sec:resultsbenchmark}. The computational results are discussed in Section~\ref{sec:resultsresults}.
Lastly, we discuss the computational times in Section~\ref{sec:resultscomptime}.

\subsection{Instance generation}\label{sec:resultsinstances}
We generate instances with different CS characteristics, such as variable capacities, densities, CS technology dispersion, and arrival rates of exogenous EVs.
Recalling that $\kappa$ is the capacity of the queue at CSs, we define $\mathcal{K}=\{1, 2, 3\}$ as the set of possible values for $\kappa$. The value of $\kappa$ is set equally for all CSs.
We consider a set $\mathcal{V}$ of three different EV types, specifically $\mathcal{V}=\{$Peugeot I-on, BMW i3s, Renault Zoe$\}$. 
Table~\ref{tab:instanceV} shows the characteristics of these three vehicles.
The consumption rate indicates the amount of energy (KWh) the EV consumes per kilometer.
\begin{table}[!htbp]
    \centering
    \footnotesize
    \begin{tabular}{cc|c|c|c}
        \multicolumn{2}{c|}{EV} & Battery size & Consumption rate  & Autonomy  \\ \cline{1-2}
        Make & Model & (KWh) & (KWh/km) & (km) \\ \hline
        Peugeot & I-on & 16 & 0.125 & 128 \\
        BMW & i3s & 42.2 & 0.165 & 255 \\
        Renault & Zoe & 52 & 0.165 & 315 \\ \hline
    \end{tabular}
    \caption{Details of the three EV types (source:~\cite{OpenEV2022})}
    \label{tab:instanceV}
\end{table}

We consider a set of three CS technologies defined as $F=\{$Slow, Normal, Fast$\}$.
As explained in Section~\ref{sec:problemdescription}, we use piecewise linear charging functions to compute the EV charging time. 
In particular, we adopt the charging time breakpoints of Slow, Normal, and Fast chargers proposed in~\cite{Montoya2017} for the Peugeot I-on and scale them for the other two considered EVs. 
Table~\ref{tab:instanceCSbp} shows the resulting charging time breakpoints for each EV and for each charger technology in $F$. 
The values in this table are the charging times $\phi_z$ (in hours), i.e., the time the vehicle should spend to charge from battery level zero to $\beta_z$.
For example, it takes 2.03 hours for a BMW to charge from zero to 95\% of its battery capacity (i.e., $\beta_z=40.1$ KWh), when charging in a Normal CS.
\begin{table}[!ht]
    \centering
    \setlength{\tabcolsep}{0.75ex}
    \footnotesize
    \begin{tabular}{c|c| c c c|c| c c c|c| c c c}
         & \multicolumn{4}{c|}{Peugeot } & \multicolumn{4}{c|}{BMW} & \multicolumn{4}{c}{Renault}\\ \cline{2-13}
        SoC (\%) & $\beta_z$ (KWh) & Slow & Normal & Fast & $\beta_z$ (KWh) & Slow & Normal & Fast & $\beta_z$ (KWh) & Slow & Normal & Fast \\ \hline
        85 & 13.6 & 1.26 & 0.62 & 0.31 & 35.9 & 3.32 & 1.64 & 0.82 & 44.2 & 4.09 & 2.02 & 1.01 \\
        95 & 15.2 & 1.54 & 0.77 & 0.39 & 40.1 & 4.06 & 2.03 & 1.03 & 49.4 & 5.00 & 2.50 & 1.27\\
        100 & 16.0 & 2.04 & 1.01 & 0.51 & 42.2 & 5.38 & 2.66 & 1.35 & 52.0& 6.63 & 3.28 & 1.66 \\ \hline 
    \end{tabular}
    \caption{Breakpoint charging time ($\phi_z$) of the piecewise linear charging functions $\Phi_i(b)$}
    \label{tab:instanceCSbp}
\end{table}

We consider a region of $350\times 700$ kilometers where the origin and destination are located in corners $o=(0, 0)$ and $d=(350, 700)$.
We assume that distances between every two locations in this region are euclidean and that the EV travels at the constant speed of 100 km/h.
We consider a set $\mathcal{D}$ of three density levels to determine the number of CSs in the region, specifically $\mathcal{D}=\{$Low, Moderate, High$\}$.
In regions with Low, Moderate, and High densities, there are 2, 6, and 10 CSs per $100\times100$ km, respectively. 
Consequently, the total number of CSs in instances with Low, Moderate, and High density is 49, 147, and 245, respectively.
The locations of the CSs are generated following a uniform distribution.
Each generated CS is assigned a charging technology from $F$, where $F$ is the set of charging technologies, according to a given probability distribution $P$.
The distribution $P$ belongs to a set $\mathcal{P}$ of six distributions, $\mathcal{P}=\{PU_1, PU_2, PU_3, PM_1, PM_2, PM_3\}$, which are described in Table~\ref{tab:instancep}. 
For example, in $PM_2$, 55\% of CSs are Slow, 42\% are Normal, and 3\% of them are Fast.
The first three distributions (i.e., $PU_1, PU_2, PU_3$) are designed to analyze the impact of unique CS technologies, whereas the latter three distributions (i.e., $PM_1, PM_2, PM_3$) are designed to analyze the impact of CS technology mixes.
These were inferred from the publication~\cite{Environment2020}, which reports the current and anticipated status of public CSs in Europe. 
Accordingly, Spain, France, and Germany are examples of $PM_1, PM_2,$ and $PM_3$, respectively.

\begin{table}[!ht]
    \centering
    \footnotesize
    \begin{tabular}{c|ccc}
        \multirow{2}{*}{$P$} & \multicolumn{3}{c}{CS technology ($f$)}\\
         & Slow & Normal & Fast \\ \hline
        $PU_1$ & 100\% & 0\% & 0\% \\
        $PU_2$ & 0\% & 100\% & 0\% \\
        $PU_3$ & 0\% & 0\% & 100\% \\ \hline
        $PM_1$ & 55\% & 42\% & 3\% \\
        $PM_2$ & 45\% & 45\% & 10\% \\
        $PM_3$ & 10\% & 65\% & 25\% \\ \hline
    \end{tabular}
    \caption{Six configurations for the portion of CS technologies in the region}
    \label{tab:instancep}
\end{table}

We estimate the expected charging time of exogenous EVs by setting it to the time required to charge an EV with an average battery size from 0\% to 85\%.
According to a database containing characteristics of about 300 EV models~\citep{OpenEV2022}, the average battery size for EVs (excluding Tesla and plugin-hybrid cars) is 46 KWh. 
By using the breakpoints provided in Table~\ref{tab:instanceCSbp} and scaling them for a battery size of 46 KWh, the expected charging time (to charge from 0\% to 85\%) for Slow, Normal, and Fast CSs are computed as 3.62, 1.78, and 0.89 hours, respectively.
Consequently, we set its service rate $\mu_i$ to 0.28, 0.56, and 1.12 vehicles per hour, for $f_i\in\{$Slow, Normal, Fast$\}.$
We compute the arrival rate of EVs at a CS $i$ as a function of the utilization rate and service rate of that CS, i.e., $\lambda_i=\mu_i\rho_i$. 
We consider a set of utilization rates $\mathcal{U}=\{$Low, Medium, High$\}$, where the utilization rate of all CSs is fixed and equal to 40\%, 65\%, and 90\%, respectively.
The arrival rate for each CS technology with each utilization rate is shown in Table~\ref{tab:instanceutil}.

\begin{table}[!htbp]
    \centering
        \begin{tabular}{c|c|c|c}
            \multirow{2}{*}{$\rho_i$} & \multicolumn{3}{c}{$f_i$} \\ \cline{2-4}
             & Slow & Normal & Fast \\\hline
            40\% & 0.11 & 0.22 & 0.45 \\
            65\% & 0.18 & 0.36 & 0.73 \\
            90\% & 0.25 & 0.50 & 1.01 \\ \hline
        \end{tabular}
        \caption{Arrival rate $\lambda_i$ (EVs per hour) at a CS with a charging technology of $f_i$ for utilization rate $\rho_i$}
        \label{tab:instanceutil}
    \end{table}

We define $\mathcal{I}=\mathcal{K}\times\mathcal{V}\times\mathcal{D}\times\mathcal{P}\times\mathcal{U}$ as the set of all instances. 
We identify each \textit{instance} in $\mathcal{I}$ by a tuple of $(\kappa, V,D,P,U)$. 
The set $\mathcal{I}$ contains a total of 486 different instances. 

As discussed in Section~\ref{sec:estpath}, we classify the behavior of the waiting time of each CS $i$ into three cases, namely A, B, and C, based on its characteristics. 
In particular, we empirically classify all 27 combinations of $\lambda_i, \mu_i,$ and $\kappa_i$ that we consider in our experiments, as shown in Table~\ref{tab:cscases}.
\begin{table}[!ht]
    \centering
    \setlength{\tabcolsep}{0.75ex}
    \footnotesize
    \begin{tabular}{c|c|ccc||ccc||ccc}
        \multicolumn{2}{c}{} & \multicolumn{3}{c}{$\kappa_i=1$} &  \multicolumn{3}{c}{$\kappa_i=2$} & \multicolumn{3}{c}{$\kappa_i=3$}  \\ \cline{3-11}
        \multicolumn{2}{c|}{} & \multicolumn{3}{c||}{$\lambda_i$} &  \multicolumn{3}{c||}{$\lambda_i$} & \multicolumn{3}{c}{$\lambda_i$}  \\ \cline{3-11}
        \multicolumn{2}{c|}{} & 0.11 & 0.18 & 0.25 & 0.22 & 0.36 & 0.50 & 0.45 & 0.73 & 1.01 \\ \hline
        \multirow{3}{*}{$\mu_i$} & 0.28 & B & B & B & B & B & B & B & A & C \\
         & 0.56 & B & B & B & B & B & B & B & A & C \\
         & 1.12 & B & B & B & B & B & B & B & A & C \\ \hline
    \end{tabular}
    \caption{Classifying CSs into cases A, B, and C with respect to the behavior of waiting times}
    \label{tab:cscases}
\end{table}

\subsection{Benchmark policy}\label{sec:resultsbenchmark}
We benchmark the performance of our methodology (i.e., \roi) against the setting where the OI information is not available.
The benchmark adopted in this case is a reoptimization method that, similar to \cite{Sweda2017}, makes decisions only when the queue is observed.
The dynamics of the problem in this benchmark can be described by using the Markov chain presented in Figure \ref{fig:mdp}, considering that the ``Indicator changes'' event is excluded from the list of events that trigger a decision epoch.
As a result, decisions are reoptimized when the EV arrives at a CS, completes waiting, or the target battery level is reached.

Since actions in the benchmark are exclusively made at the origin or CSs, $N_k=N$ and $c\in N$.
When the EV leaves the origin, it has no information regarding CSs other than their asymptotic behavior. 
Therefore, the steady-state waiting time $W_j$ for each CS $j\in \mathcal{C}$ is used when computing $\Tilde{\tau}_{ij}$. 
Specifically, $\Tilde{\tau}_{ij}= t_{ij} + W_j,~\forall i, j\in N$.
Once an event is triggered, $c$ refers to a CS (i.e., $c\in\mathcal{C}$) with an observed queue length of $\xi_{c}$.
Therefore, the benchmark method updates $\Tilde{\tau}_{cc}=\xi_{c}/\mu_{c}$ when determining $\Tilde{\tau}_{ij}$.
Then, similar to the method proposed in Section~\ref{sec:solution}, the benchmark method adopts the \texttt{frvcpy} package to reoptimize the route by incorporating the updated $\Tilde{\tau}_{ij}$ matrix.
We refer to this benchmark as the Reoptimization method with the Steady-State information (\rss{}).

\subsection{Results and discussion}\label{sec:resultsresults}
To analyze the performance of our method from different angles, we consider six experimental settings.
In particular, we distinguish the instance $I_B=\{1\}\times\{\textrm{BMW}\}\times\{\textrm{Moderate}\}\times\{PU_2\}\times\{\textrm{Medium}\}$ as a baseline instance and vary it with respect to different instance identifiers at each experimental setting.  
The first experimental setting aims at analyzing the performance of \roi{} and \rss{} with respect to the CS technologies. 
To this end, we solve instances $\mathcal{I}_{PU}=\{1\}\times\{\textrm{BMW}\}\times\{\textrm{Moderate}\}\times\{PU_1, PU_2, PU_3\}\times\{\textrm{Medium}\}$. 
In the second experimental setting, we explore more realistic mixes of CS technologies by solving instances $\mathcal{I}_{PM}=\{1\}\times\{\textrm{BMW}\}\times\{\textrm{Moderate}\}\times\{PM_1, PM_2, PM_3\}\times\{\textrm{Medium}\}$. 
In the third experimental setting, we explore the different values of $\kappa$, where we consider instances  $\mathcal{I}_{\kappa}=\mathcal{K}\times\{\textrm{BMW}\}\times\{\textrm{Moderate}\}\times\{PU_2\}\times\{\textrm{Medium}\}$.
We investigate the performance of our method considering the three vehicle types in the forth experimental setting, i.e., $\mathcal{I}_{V}=\{1\}\times\mathcal{V}\times\{\textrm{Moderate}\}\times\{PU_2\}\times\{\textrm{Medium}\}$. 
In the fifth experimental setting, we examine the effect of the different density levels by solving instances $\mathcal{I}_{D}=\{1\}\times\{\textrm{BMW}\}\times\mathcal{D}\times\{PU_2\}\times\{\textrm{Medium}\}$.
Finally, in the last experimental setting, we solve instances $\mathcal{I}_{U}=\{1\}\times\{\textrm{BMW}\}\times\{\textrm{Moderate}\}\times\{PU_2\}\times\mathcal{U}$ to investigate the impact of different utilization rates.
A summary of the experimental settings is reported in Table~\ref{tab:allinstances}.

We recall that CS locations are assumed to be uniformly distributed for a given density level. Therefore, for consistency reasons, we sample one realization of CS locations for each of the three density levels and use it whenever needed.
Specifically, a sample set of CS locations for the Moderate density level is generated. This sample is then modified by uniformly sampling additional CS locations to create the sample set for the High density level. To create the Low density sample, CSs are randomly removed from the Moderate sample. 
In the second experimental setting, where $PM_1, PM_2$, and $PM_3$ refer to distribution functions for CS technologies, we sample one realization for each.
Following a similar procedure as the one used to realize CS locations for different density levels, we sample one CS technologies assignment for $PM_2$ under its distribution function and then modify it to create a sample  for $PM_1$ and $PM_3$. 
Given that the portion of Slow, Normal, and Fast CSs in $PM_2$ are 45\%, 45\%, and 10\%, respectively, to generate the sample for $PM_1$, we randomly select 3\% of Normal and 7\% of Fast CSs in $PM_2$ and convert them to Slow CSs.
Similarly, in order to generate the sample for $PM_3$, we first select 35\% of CSs with Slow charging technology, we then convert 20\% of them to Normal and 15\% to Fast CSs.
Lastly, we test the performance of \roi{} and \rss{} on each instance for 15 daily scenarios.
Each daily scenario corresponds to sampling a set of inter-arrival and service times of exogenous EVs at each CS according to Table~\ref{tab:instanceutil}.
We note that in nine out of all 18 instances that we consider in the experimental settings, we have $D=$ Moderate and $P=PU_2$, meaning that the set of CS locations and their technologies are identical. Therefore, 
the same daily scenarios are used for them.
Thus, a total of 10 samples, each containing 15 days, are used in the experiments.
We note that the same daily scenarios are used for \roi{} and \rss{}.
\begin{table}[!ht]
    \centering
    \setlength{\tabcolsep}{0.75ex}
    \footnotesize
    \begin{tabular}{C{1.8cm}|c|ccccc}
        Experimental setting & Instances & $\kappa$ & $V$ & $D$ & $P$ & $U$ \\ \hline
         & & & & & $PU_1$ & \\
        1 & $\mathcal{I}_{PU}$ & 1 & BMW & Moderate & $PU_2$ & Medium \\
         & & & & & $PU_3$ & \\ \hline
         & & & & & $PM_1$ & \\
        2 & $\mathcal{I}_{PM}$ & 1 & BMW & Moderate & $PM_2$ & Medium \\
         & & & & & $PM_3$ & \\ \hline
         & & 1 & & & & \\
        3 & $\mathcal{I}_{\kappa}$ & 2 & BMW & Moderate & $PU_2$ & Medium \\
         & & 3 & & & & \\ \hline
         & & & Peugeot & & & \\
        4 & $\mathcal{I}_{V}$ & 1 & BMW & Moderate & $PU_2$ & Medium \\
         & & & Renault & & & \\ \hline
         & & & & Low & & \\
        5 & $\mathcal{I}_{D}$ & 1 & BMW & Moderate & $PU_2$ & Medium \\
         & & & & High & & \\ \hline 
         & & & & & & Low \\
        6 & $\mathcal{I}_{U}$ & 1 & BMW & Moderate & $PU_2$ & Medium \\
         & & & & & & High \\ \hline 
    \end{tabular}
    \caption{Experimental settings and the underlying instances}
    \label{tab:allinstances}
\end{table}

The results for each of the previously described experimental settings are presented in two sets of tables.
In the first set of tables (i.e., Tables~\ref{tab:resp012a},~\ref{tab:resp345a},~\ref{tab:reskappaa},~\ref{tab:resvehiclea},~\ref{tab:resdensitya}, and~\ref{tab:resutilizationa}), we report the waiting time (Wait), charging time (Charge), driving time (Drive), and total trip duration (Total) averaged over the 15 daily scenarios for \roi{} and \rss{}.
We note that the duration is reported in HH:MM format.
The first column of the ``Improvement" section  reports the average waiting time improvement ($\Delta$Wait), computed as Wait in \roi{} minus Wait in \rss{}. 
This section also reports the average percentage improvements in waiting time (Wait), charging time (Charge), driving time (Drive), and total trip duration (Total).
We compute the percentage improvement of \roi{} over \rss{} with respect to the factor $X$ as follows:
$\displaystyle{({X \textrm{(\roi{})} - X \textrm{(\rss{})}})/{X \textrm{(\rss{})}}\times100}$.
In the second set of tables (i.e., Tables~\ref{tab:resp012b},~\ref{tab:resp345b},~\ref{tab:reskappab},~\ref{tab:resvehicleb},~\ref{tab:resdensityb}, and~\ref{tab:resutilizationb}), we report the average number of times per day the EV: performs charging (\chr{}), deviates towards another CS while driving (\dd{}), and deviates towards another CS while waiting (or just arrived) at a CS (\dc{}) for \roi{} and \rss{}. 

In the first experimental setting, we compare the performance of \roi{} with \rss{} on instances with all Slow ($PU_1$), all Normal ($PU_2$), and all Fast ($PU_3$) CSs, i.e., $\mathcal{I}_{PU}$.
Tables~\ref{tab:resp012a} and~\ref{tab:resp012b} summarize these results.
Table~\ref{tab:resp012a} demonstrates that the charging time, and total trip time, for both methods, decrease significantly as the charging technology of the CSs becomes faster.
Comparing the two methods shows that \roi{} decreases the waiting time and the total trip duration by an average of 89.8\% and 10.8\%, respectively. 
According to this table, the advantage of \roi{} over \rss{} is more pronounced on instances with Slow CSs.
Furthermore, whereas the average waiting time in \rss{} tends to be considerably high when the charging technology is Slow or Normal, our method maintains the waiting time at a relatively consistent low levels.
\begin{table}[!hbtp]
\centering
\setlength{\tabcolsep}{0.75ex}
\footnotesize
\begin{tabular}{c|cccc|cccc|ccccc}
 & \multicolumn{4}{c}{\rss} & \multicolumn{4}{c}{\roi} & \multicolumn{5}{c}{Improvement} \\
$P$ & Wait & Charge & Drive & Total & Wait & Charge & Drive & Total  & $\Delta$Wait & Wait & Charge & Drive & Total \\ \hline
$PU_1$   & 3:17    & 8:15     & 8:06    & 19:38 & 0:00    & 8:29     & 8:01    & 16:30 & -3:17 & -99.5\% & 2.9\%    & -1.0\%  & -15.8\% \\
$PU_2$ & 1:48    & 4:03     & 8:04    & 13:55 & 0:08    & 4:09     & 7:58    & 12:15 & -1:40 & -92.0\% & 2.6\%    & -1.2\%  & -11.9\% \\
$PU_3$   & 0:53    & 1:59     & 7:50    & 10:42 & 0:11    & 2:03     & 7:57    & 10:11 & -0:42 & -77.7\% & 3.3\%    & 1.6\%   & -4.7\%  \\ \hline
\textbf{Avg.}       & \textbf{1:59}    & \textbf{4:46}     & \textbf{8:00}    & \textbf{14:46} & \textbf{0:07}    & \textbf{4:54}     & \textbf{7:59}    & \textbf{13:00} & \textbf{-1:52} & \textbf{-89.8\%} & \textbf{2.9\%}    & \textbf{-0.2\%}  & \textbf{-10.8\%} \\ \hline
\end{tabular}
\caption{Results for experimental setting 1 ($\mathcal{I}_{PU}$)}
\label{tab:resp012a}
\end{table}

According to Table~\ref{tab:resp012b}, the EV in \rss{} leaves a visited CS without charging an average of 0.3 times per day, whereas this value is zero in \roi{}. This result can be explained as follows. 
Recalling that \rss{} makes routing decisions purely relying on steady-state waiting times, the EV in this method may observe the queue length of a CS for which the corresponding expected waiting time is considerably higher than anticipated.
As a result, the EV is re-routed towards a different CS. 
On the contrary, \roi{} takes advantage of the OI updates to modify the routing decision while driving, and this results in a lower probability of visiting a busy CS with longer-than-anticipated expected waiting times upon arrival.
We also observe that the number of deviations while driving in \roi{} increases with the speed of the charging technology. 
This may be related to the way our instances are generated. In fact, the utilization rate of each CS is kept constant through $PU_1$, $PU_2$ and $PU_3$. 
Therefore, instances with faster technologies also have higher arrival rates and more frequent OI updates, which leads to increasing the number of routing decisions in \roi{}. 
We recall that deviation while driving is not permitted in \rss{}, consequently, column \dd{} in \rss{} always reports zero. 

\begin{table}[!hbtp]
\centering
\setlength{\tabcolsep}{0.75ex}
\footnotesize
\begin{tabular}{c|ccc|ccc}
 & \multicolumn{3}{c}{\rss} & \multicolumn{3}{c}{\roi} \\
$P$ & \chr{} & \dd{} & \dc{} & \chr{} & \dd{} & \dc{} \\ \hline
$PU_1$ & 3.0          & 0.0          & 0.5          & 3.1          & 2.1          & 0.0 \\
$PU_2$ & 3.0          & 0.0          & 0.5          & 3.1          & 3.3          & 0.0 \\
$PU_3$ & 3.0          & 0.0          & 0.0          & 3.1          & 4.5          & 0.0 \\ \hline
\textbf{Avg.} & \textbf{3.0} & \textbf{0.0} & \textbf{0.3} & \textbf{3.1} & \textbf{3.3} & \textbf{0.0}  \\ \hline
\end{tabular}
\caption{Performance indicators for experimental setting  1 ($\mathcal{I}_{PU}$)}
\label{tab:resp012b}
\end{table}

In the second experimental setting, we study the performance of \roi{} on instances in $\mathcal{I}_{PM}$ (i.e., instances with different mixes of CS technologies). 
Results of this experimental setting are reported in Tables~\ref{tab:resp345a} and~\ref{tab:resp345b}.
Recall that $PM_1$ corresponds to instances with a lower percentage of Fast chargers, where the majority of CSs are Slow, whereas $PM_3$ corresponds to instances with a higher percentage of Fast chargers, where the majority of CSs are Normal.
Comparing the two methods in Table~\ref{tab:resp345a} shows that \roi{} decreases the waiting time and total trip duration by 23.7\% and 1.4\%, respectively. 
Furthermore, when comparing $PM_1$ to $PM_2$ and $PM_2$ to $PM_3$, as the proportion of CSs with a Fast charging technology increases, the charging, driving, and total trip time decrease in both methods.
A comparison between this experimental setting and the previous one shows that the presence of Fast CSs implies a reduction in the total trip duration for both methods. However, the advantage of \roi{} over \rss{} is considerably less pronounced in the second experimental setting. A detailed analysis of the results, shows that, as expected, in \rss{} the EV charges at Fast CSs, whenever possible. Similar solutions are found by \roi{} because selecting a Fast CS is a good option even when the estimated expected waiting time is relatively high. 
In other words, Fast CSs are so fast that whenever there is a mix of Slow, Normal and Fast CSs, charging at Fast CSs, even if busy, is on average much more convenient than charging at an available Slow or Normal CS. This reasoning is confirmed by the fact that the average waiting time for \roi{} in Table~\ref{tab:resp345a} is higher than what was observed for the previous experiment in Tables~\ref{tab:resp012a}.
For the same reasons, the average charging time in  presence of Fast CSs is considerably shorter than what was observed in Tables~\ref{tab:resp012a} for $PU_1$ to $PU_2$, which do not have Fast CSs.
\begin{table}[!hbtp]
\centering
\setlength{\tabcolsep}{0.75ex}
\footnotesize
\begin{tabular}{c|cccc|cccc|ccccc}
 & \multicolumn{4}{c}{\rss} & \multicolumn{4}{c}{\roi} & \multicolumn{5}{c}{Improvement} \\
$P$ & Wait & Charge & Drive & Total & Wait & Charge & Drive & Total  & $\Delta$Wait & Wait & Charge & Drive & Total \\ \hline
$PM_1$ & 1:03 & 2:47 & 7:57 & 11:48 & 0:29 & 2:52 & 8:02 & 11:24 & -0:33 & -52.7\% & 2.5\% & 1.0\% & -3.4\%          \\
$PM_2$ & 1:16 & 2:07 & 7:56 & 11:19 & 1:12 & 2:08 & 7:57 & 11:18 & -0:03 & -4.6\% & 1.4\% & 0.1\% & -0.2\%          \\
$PM_3$ & 1:04 & 1:54 & 7:55 & 10:54 & 0:55 & 2:01 & 7:54 & 10:51 & -0:08 & -13.8\% & 5.5\% & -0.2\% & -0.5\%          \\ \hline
\textbf{Avg.} & \textbf{1:07} & \textbf{2:16} & \textbf{7:56} & \textbf{11:21} & \textbf{0:52} & \textbf{2:20} & \textbf{7:58} & \textbf{11:11} & \textbf{-0:15} & \textbf{-23.7\%} & \textbf{3.1\%} & \textbf{0.3\%} & \textbf{-1.4\%} \\ \hline
\end{tabular}
\caption{Results for experimental setting 2 ($\mathcal{I}_{PM}$)}
\label{tab:resp345a}
\end{table}

Results in Table~\ref{tab:resp345b} indicate that both methods have the most deviations in $PM_3$. As previously explained, the EV generally prefers visiting Fast CSs. Therefore, when there is a limited number of Fast CSs (such as in $PM_1$ or $PM_2$), the options to deviate are reduced. 
\begin{table}[!hbtp]
\centering
\setlength{\tabcolsep}{0.75ex}
\footnotesize
\begin{tabular}{c|ccc|ccc}
 & \multicolumn{3}{c}{\rss} & \multicolumn{3}{c}{\roi} \\
$P$ & \chr{} & \dd{} & \dc{} & \chr{} & \dd{} & \dc{}  \\ \hline
$PM_1$ & 3.0          & 0.0          & 0.0            & 3.2          & 0.7          & 0.1            \\
$PM_2$ & 3.0          & 0.0          & 0.0            & 3.1          & 0.1          & 0.0             \\
$PM_3$ & 3.0          & 0.0          & 0.7            & 3.1          & 1.7          & 0.1             \\ \hline
\textbf{Avg.} & \textbf{3.0} & \textbf{0.0} & \textbf{0.2} & \textbf{3.1} & \textbf{0.8} & \textbf{0.1}  \\ \hline
\end{tabular}
\caption{Performance indicators for experimental setting 2 ($\mathcal{I}_{PM}$)}
\label{tab:resp345b}
\end{table}

Tables~\ref{tab:reskappaa} and~\ref{tab:reskappab} show the results for $\mathcal{I}_{\kappa}$, which  correspond to the third experimental setting.
As shown in Table~\ref{tab:reskappaa}, \roi{} reduces the overall trip duration by an average of 18.5\%. 
This is mainly due to an average reduction of 3:17 in waiting time.
Table~\ref{tab:reskappaa} shows that as $\kappa$ becomes larger, the waiting time and the total trip duration tend to increase in both methods. 
This finding can be explained by the fact that, in all types of queues that we consider, an increase in the number of EVs allowed to wait in a CS corresponds to a longer average waiting time, thus leading to longer total trip durations on average. 
However, compared to \roi{}, \rss{} shows higher increases in waiting time as $\kappa$ increases. 
In particular, as shown in column $\Delta$Wait, the absolute gap between \roi{} and \rss{} increases with $\kappa$.
Therefore, we infer that the effectiveness of our method is amplified in instances with larger $\kappa$.
According to Table~\ref{tab:reskappaa}, the average charging time in \roi{} is 7.3\% longer than that in \rss{}. This variation may be linked to the higher number of charging sessions (\chr{}) that the EV performs in \roi{}, as shown in Table \ref{tab:reskappab}. In particular, since the difference between the driving time in both methods is not significant (i.e., the EV in both methods have consumed and charged about the same amount of energy), we infer that the EV in \roi{} has arrived at CSs with a higher battery level, which led to slower charging rates.

\begin{table}[!ht]
\centering
\setlength{\tabcolsep}{0.75ex}
\footnotesize
\begin{tabular}{c|cccc|cccc|ccccc}
 & \multicolumn{4}{c}{\rss} & \multicolumn{4}{c}{\roi} & \multicolumn{5}{c}{Improvement} \\
$\kappa$ & Wait & Charge & Drive & Total & Wait & Charge & Drive & Total  & $\Delta$Wait & Wait & Charge & Drive & Total \\ \hline
1 & 1:48 & 4:03 & 8:04 & 13:55 & 0:08 & 4:09 & 7:58 & 12:15 & -1:40 & -92.0\% & 2.6\% & -1.2\% & -11.9\%          \\
2 & 3:58 & 4:00 & 7:54 & 15:52 & 0:06 & 4:23 & 8:06 & 12:35 & -3:52 & -97.5\% & 9.9\% & 2.5\% & -20.6\%          \\
3 & 4:37 & 4:01 & 8:00 & 16:38 & 0:18 & 4:24 & 8:08 & 12:50 & -4:19 & -93.3\% & 9.3\% & 1.7\% & -22.9\%          \\ \hline
\textbf{Avg.} & \textbf{3:28} & \textbf{4:01} & \textbf{7:59} & \textbf{15:29} & \textbf{0:11} & \textbf{4:19} & \textbf{8:04} & \textbf{12:34} & \textbf{-3:17} & \textbf{-94.3\%} & \textbf{7.3\%} & \textbf{1.0\%} & \textbf{-18.5\%} \\ \hline
\end{tabular}
\caption{Results for experimental setting 3 ($\mathcal{I}_{\kappa}$)}
\label{tab:reskappaa}
\end{table}

Table~\ref{tab:reskappab} shows that \roi{} increases the daily average number of charging sessions (\chr{}) from 3.0 in \rss{} to 3.4. Additionally, while the number of charging sessions per day remains constant at 3.0 for $\kappa=1, 2, 3$ in \rss{}, it increases from 3.1 to 3.7 in \roi{} as $\kappa$ increases from 1 to 3. One possible explanation for this finding is that higher values of $\kappa$ increase the average queue lengths.
As all CSs have the same technology, \rss{} treats them all similarly, rendering $\kappa$ change ineffective in its decision-making process. Specifically, \rss{} yields the minimum number of charging sessions per route, which is 3.0 in all cases. In contrast, as our method computes waiting times in a more refined fashion, it may choose to perform a larger number of charging sessions if the corresponding total estimated waiting time is shorter. 
Similar to the first experimental setting, this table indicates that in \roi{} the EV rarely leaves a visited CS (\dc{} is almost zero on average), while we observe a considerably higher number of deviations while driving (on average of 3.8). This implies that, even though the expected queue length of CSs in instances with a higher $\kappa$ is larger, the waiting time estimations made en route align with the actual statuses of the queues.
\begin{table}[!ht]
\centering
\setlength{\tabcolsep}{0.75ex}
\footnotesize
\begin{tabular}{c|ccc|ccc}
 & \multicolumn{3}{c}{\rss} & \multicolumn{3}{c}{\roi} \\
$\kappa$ & \chr{} & \dd{} & \dc{} & \chr{} & \dd{} & \dc{} \\ \hline
1 & 3.0 & 0.0 & 0.5 & 3.1 & 3.3 & 0.0 \\
2 & 3.0 & 0.0 & 0.1 & 3.5 & 4.2 & 0.0 \\
3 & 3.0 & 0.0 & 0.3 & 3.7 & 3.8 & 0.1 \\ \hline
\textbf{Avg.} & \textbf{3.0} & \textbf{0.0} & \textbf{0.3} & \textbf{3.4} & \textbf{3.8} & \textbf{0.0}  \\ \hline
\end{tabular}
\caption{Performance indicators for experimental setting 3 ($\mathcal{I}_{\kappa}$)}
\label{tab:reskappab}
\end{table}

The results of the fourth experimental setting $\mathcal{I}_V$ are reported in Tables~\ref{tab:resvehiclea} and~\ref{tab:resvehicleb}.
We recall that Peugeot and Renault have the smallest and largest autonomy, respectively. 
According to Table~\ref{tab:resvehiclea}, as the EV autonomy increases (e.g., BMW to Renault), the waiting time and total trip duration decrease in both methods.
This is because EVs with greater autonomy require fewer charging sessions and are therefore less likely to experience waiting at CSs. 
Results also show that the improvement of \roi{} over \rss{} decreases as the EV autonomy increases, confirming that our method is more effective when EVs require more frequent charging.
We note that the driving time for both methods decreases with the increase in the EV autonomy. This may be explained by the fact that greater autonomy entails less deviations to CSs.
\begin{table}[!ht]
\centering
\setlength{\tabcolsep}{0.75ex}
\footnotesize
\begin{tabular}{c|cccc|cccc|ccccc}
 & \multicolumn{4}{c}{\rss} & \multicolumn{4}{c}{\roi} & \multicolumn{5}{c}{Improvement} \\
$V$ & Wait & Charge & Drive & Total & Wait & Charge & Drive & Total  & $\Delta$Wait & Wait & Charge & Drive & Total \\ \hline
Peugeot & 6:14          & 3:49          & 8:08          & 18:11          & 0:38          & 4:13          & 8:09          & 13:00          & -5:36          & -89.7\%          & 10.6\%         & 0.3\%           & -28.4\%          \\
BMW           & 1:48          & 4:03          & 8:04          & 13:55          & 0:08          & 4:09          & 7:58          & 12:15          & -1:40          & -92.0\%          & 2.6\%          & -1.2\%          & -11.9\%          \\
Renault       & 0:25          & 3:29          & 7:54          & 11:48          & 0:06          & 3:35          & 7:53          & 11:34          & -0:19          & -73.0\%          & 2.8\%          & -0.1\%          & -1.9\%           \\ \hline
\textbf{Avg.} & \textbf{2:49} & \textbf{3:47} & \textbf{8:02} & \textbf{14:38} & \textbf{0:18} & \textbf{3:59} & \textbf{8:00} & \textbf{12:17} & \textbf{-2:31} & \textbf{-84.9\%} & \textbf{5.4\%} & \textbf{-0.3\%} & \textbf{-14.1\%}  \\ \hline
\end{tabular}
\caption{Results for experimental setting 4 ($\mathcal{I}_{V}$)}
\label{tab:resvehiclea}
\end{table}

Similar to the previous experimental settings, Table~\ref{tab:resvehicleb} shows the average \dc{} is reduced from 0.9 in \rss{} to 0.1 in \roi{}, which implies that \roi{} properly adjusts the EV while driving so that it does not need to leave the CS after observing the actual queue length. Additionally, as autonomy increases, the average \dd{} increases. A possible explanation for this finding is that as autonomy increases, the EV will have a larger set of potential CSs towards which it may  deviate to.
\begin{table}[!ht]
\centering
\setlength{\tabcolsep}{0.75ex}
\footnotesize
\begin{tabular}{c|ccc|ccc}
 & \multicolumn{3}{c}{\rss} & \multicolumn{3}{c}{\roi} \\
$V$ & \chr{} & \dd{} & \dc{} & \chr{} & \dd{} & \dc{} \\ \hline
Peugeot       & 7.0          & 0.0          & 1.2          & 7.4          & 2.5          & 0.2 \\
BMW           & 3.0          & 0.0          & 0.5          & 3.1          & 3.3          & 0.0 \\
Renault       & 2.0          & 0.0          & 1.1          & 2.0          & 3.7          & 0.0 \\ \hline
\textbf{Avg.} & \textbf{4.0} & \textbf{0.0} & \textbf{0.9} & \textbf{4.2} & \textbf{3.2} & \textbf{0.1}  \\ \hline
\end{tabular}
\caption{Performance indicators for experimental setting 4 ($\mathcal{I}_{V}$)}
\label{tab:resvehicleb}
\end{table}

Tables \ref{tab:resdensitya} and \ref{tab:resdensityb} present the results of the fifth experimental setting, which examines the performance of the proposed method on instances $\mathcal{I}_D$. 
According to Table \ref{tab:resdensitya}, \roi{} outperforms \rss{} by reducing the waiting time and total trip duration by an average of 95.4\% and 10.4\%, respectively. 
The waiting time in both methods remains relatively constant as $D$ changes from Low to High, with no direct effect on the total trip time. This can be explained by the fact that the expected waiting time at each CS remains the same at all levels of density, since CS technologies, utilization rates, and $\kappa$ are the same. On the other hand, driving time in both methods decreases as density increases. In particular, due to the larger number of CSs in instances with higher density, the EV performs shorter deviations to reach a CS.
\begin{table}[!ht]
\centering
\setlength{\tabcolsep}{0.75ex}
\footnotesize
\begin{tabular}{c|cccc|cccc|ccccc}
 & \multicolumn{4}{c}{\rss} & \multicolumn{4}{c}{\roi} & \multicolumn{5}{c}{Improvement} \\
$D$ & Wait & Charge & Drive & Total & Wait & Charge & Drive & Total  & $\Delta$Wait & Wait & Charge & Drive & Total \\ \hline
Low    & 1:36          & 4:04          & 8:08          & 13:48          & 0:05          & 4:07          & 7:58          & 12:10          & -1:31          & -94.1\%          & 1.4\%          & -2.0\%          & -11.7\%          \\
Moderate        & 1:48          & 4:03          & 8:04          & 13:55          & 0:08          & 4:09          & 7:58          & 12:15          & -1:40          & -92.0\%          & 2.6\%          & -1.2\%          & -11.9\%          \\
High          & 1:26          & 3:49          & 7:52          & 13:07          & 0:00          & 4:11          & 7:56          & 12:07          & -1:26          & -100.0\%         & 9.6\%          & 0.7\%           & -7.7\%           \\ \hline
\textbf{Avg.} & \textbf{1:37} & \textbf{3:59} & \textbf{8:01} & \textbf{13:37} & \textbf{0:04} & \textbf{4:09} & \textbf{7:57} & \textbf{12:10} & \textbf{-1:33} & \textbf{-95.4\%} & \textbf{4.6\%} & \textbf{-0.8\%} & \textbf{-10.4\%} \\ \hline
\end{tabular}
\caption{Results for experimental setting 5 ($\mathcal{I}_{D}$)}
\label{tab:resdensitya}
\end{table}

According to Table~\ref{tab:resdensityb}, the EV in \roi{} performs an average of 3.3 charging sessions, while it performs an average of 3.0 charging sessions in \rss{}. 
As in the third experimental setting, the increase in \chr{} is associated with the tendency of \roi{} to yield routes with a slightly higher number of CS visits, as this leads to a lower total expected waiting time.

The table also shows that for \rss{} the average number of deviations while waiting (\dc{}) in the High density case is double that of the Low density.
A possible explanation for this is that the EV has a larger set of potential CSs to deviate to in higher density environments. Additionally, as $D$ increases, CSs become closer to each other, resulting in smaller deviation costs for the EV. This reasoning also explains why \dd{} in \roi{} increases relatively with $D$. 
\begin{table}[!ht]
\centering
\setlength{\tabcolsep}{0.75ex}
\footnotesize
\begin{tabular}{c|ccc|ccc}
 & \multicolumn{3}{c}{\rss} & \multicolumn{3}{c}{\roi} \\
$D$ & \chr{} & \dd{} & \dc{} & \chr{} & \dd{} & \dc{} \\ \hline
Low           & 3.0          & 0.0          & 0.6         & 3.3          & 3.3          & 0.2            \\
Moderate        & 3.0          & 0.0          & 0.5          & 3.1          & 3.3          & 0.0               \\
High          & 3.0          & 0.0          & 1.0            & 3.5          & 4.2          & 0.0           \\ \hline
\textbf{Avg.} & \textbf{3.0} & \textbf{0.0} & \textbf{0.7} & \textbf{3.3} & \textbf{3.6} & \textbf{0.1}  \\ \hline
\end{tabular}
\caption{Performance indicators for experimental setting 5 ($\mathcal{I}_{D}$)}
\label{tab:resdensityb}
\end{table}

Tables~\ref{tab:resutilizationa} and~\ref{tab:resutilizationb} present the results of the last experimental setting, which examines different utilization rates (i.e., instances $\mathcal{I}_U$). 
Table~\ref{tab:resutilizationa} shows that \roi{} decreases the waiting time and total trip duration by an average of 87.7\% and 10.1\%, respectively. 
As the utilization rate increases from Low to High, both methods show an increase in waiting time and total trip duration. 
An explanation for this is that instances with higher utilization rates involve a greater number of exogenous vehicles, which increases the likelihood of the EV arriving at busy CSs, leading to potentially longer waiting times. 
According to Table~\ref{tab:resutilizationb}, increasing the utilization rate leads to more frequent deviations of the EV. This may be explained by the fact that the OI information is updated more frequently with higher arrival rates of exogenous EVs.
\begin{table}[!ht]
\centering
\setlength{\tabcolsep}{0.75ex}
\footnotesize
\begin{tabular}{c|cccc|cccc|ccccc}
 & \multicolumn{4}{c}{\rss} & \multicolumn{4}{c}{\roi} & \multicolumn{5}{c}{Improvement} \\
$U$ & Wait & Charge & Drive & Total & Wait & Charge & Drive & Total  & $\Delta$Wait & Wait & Charge & Drive & Total \\ \hline
Low           & 0:59          & 4:01          & 7:58          & 12:58          & 0:03          & 4:04          & 7:54          & 12:01          & -0:56          & -93.5\%          & 1.4\%          & -0.7\%          & -7.1\%           \\
Medium      & 1:48          & 4:03          & 8:04          & 13:55          & 0:08          & 4:09          & 7:58          & 12:15          & -1:40          & -92.0\%          & 2.6\%          & -1.2\%          & -11.9\%          \\
High          & 1:53          & 4:05          & 8:12          & 14:10          & 0:25          & 4:11          & 7:58          & 12:34          & -1:28          & -77.6\%          & 2.3\%          & -2.8\%          & -11.3\%          \\ \hline
\textbf{Avg.} & \textbf{1:34} & \textbf{4:03} & \textbf{8:04} & \textbf{13:41} & \textbf{0:12} & \textbf{4:08} & \textbf{7:57} & \textbf{12:17} & \textbf{-1:22} & \textbf{-87.7\%} & \textbf{2.1\%} & \textbf{-1.6\%} & \textbf{-10.1\%} \\ \hline
\end{tabular}
\caption{Results for experimental setting 6 ($\mathcal{I}_{U}$)}
\label{tab:resutilizationa}
\end{table}

\begin{table}[!ht]
\centering
\setlength{\tabcolsep}{0.75ex}
\footnotesize
\begin{tabular}{c|ccc|ccc}
 & \multicolumn{3}{c}{\rss} & \multicolumn{3}{c}{\roi} \\
$U$ & \chr{} & \dd{} & \dc{} & \chr{} & \dd{} & \dc{} \\ \hline
Low           & 3.0          & 0.0          & 0.3          & 3.0          & 1.8          & 0.0           \\
Medium      & 3.0          & 0.0          & 0.5          & 3.1          & 3.3          & 0.0           \\
High          & 3.0          & 0.0          & 0.7          & 3.3          & 4.4          & 0.0           \\ \hline
\textbf{Avg.} & \textbf{3.0} & \textbf{0.0} & \textbf{0.5} & \textbf{3.1} & \textbf{3.2} & \textbf{0.0}  \\ \hline
\end{tabular}
\caption{Performance indicators for experimental setting 6 ($\mathcal{I}_{U}$)}
\label{tab:resutilizationb}
\end{table}

\subsection{Computation time}\label{sec:resultscomptime}
In this study, all procedures were implemented in Python and executed on 3.6 GHz Intel Xeon processors with 64 GB of RAM. 
The computation time of the proposed reoptimization algorithm is mainly influenced by two aspects: the number of times the reoptimization process is repeated (i.e., the number of decision epochs), and the duration of each reoptimization problem. 
The number of decision epochs that occur during the trip from the origin to the destination heavily depends on the parameters of the problem, including $\kappa, D, P,$ and $U$. When $\kappa$ is smaller or $D, P, U$ are at higher levels, the total number of decision epochs increases.  
In our experiments, we found that the average number of decision epochs ranges from 200 to 800. 
In practice, it is possible to control the number of decision epochs by grouping them into fixed time intervals. The reoptimization is then performed once per interval, if at least one event occurred in that interval.

Recall that each reoptimization problem in \roi{} includes both estimating the waiting time and solving the deterministic EVSPP using the \texttt{frvcpy} package.
We observed that the computation time of each reoptimization problem depends on the density of CSs in the region.
Accordingly, the average computation time for each decision epoch (including the waiting time estimation and the \texttt{frvcpy} package execution) in instances with Low, Moderate, and High densities is 0.6, 1.3, and 3.8 seconds, respectively.

\section{Conclusions}\label{sec:conclusion}
In this paper, we introduced the EVSPP-OI. This problem deals with an EV having to travel from an origin to a destination in the shortest amount of time. The EV may charge its battery at public CSs, with uncertain waiting times. However, the status of each CS is updated in real-time via binary occupancy indicator (i.e., OI) information signaling if a CS is busy or not. We examine the added value of incorporating this OI information in making routing decisions.

We formalized the \evspp{} as an MDP and solved the corresponding problem via a reoptimization algorithm, which we denote by \roi{}. In particular, at each system update (i.e., when either the EV observes a queue or the OI changes), the algorithm generates a sequence of visits to CSs and establishes their associated amount of charge. The EV follows this sequence until the next system update, at which point the procedure is repeated. The proposed method is made up of three main components. The first component is a preprocessing phase that simulates the dynamics of the queues. In particular, for each CS, this component estimates the expected waiting time as a function of the predicted arrival time  and the time elapsed since the last  relevant system update. Given the possibility of the EV making stops at intermediate CSs for charging and potentially waiting, accurately predicting the  arrival time at a specific CS presents a significant challenge. Therefore, the second component is an efficient heuristic  that estimates the arrival time of the EV at each CS. Given the estimation performed in these two components, we define  a compatible  deterministic EVSPP, which we solve with an existing algorithm.

We conducted a comprehensive computational study comprised of six experimental settings. Specifically, 
we evaluated the performance of \roi{} by comparing it with an alternative algorithm that does not use  OI information, which we denote by \rss{}.
The computational results showed that \roi{} outperformed \rss{} by 1.4\%-18.5\% and 23.7\%-95.4\% in terms of the total trip duration and waiting times, respectively. Overall \roi{} was more effective in instances with slower CSs, higher values of $\kappa$, and when more charging sessions per trip were required. 
Furthermore, we observed a number of key insights from our computational results. First, contrary to the \rss{} solutions,  in the \roi{} solutions the EV almost never leaves a CS. This is due to the fact that \roi{} efficiently reroutes the EV upon receiving OI updates. Second, compared  to \rss{}, in several cases the \roi{} solutions entailed  performing  a larger number of charging sessions. However, this is counterbalanced by a significant reduction in waiting time. Third,  the number of deviations from the planned route in \roi{} increased in instances with faster CS technology, higher EV autonomy, greater CS density, and higher utilization rates. 

Further research efforts could explore mechanisms to aggregate the CSs into clusters. This may accelerate the algorithm on larger networks. 
Additionally, while our methodology allows having several independent CSs (with single servers) at the same location,  expanding it to deal with multiple server queues is a promising research direction.

\bibliographystyle{elsarticle-harv} 
\bibliography{references.bib}
\end{document}